\documentclass{article}
\usepackage[utf8]{inputenc}
\usepackage[margin=1.0in]{geometry}
\usepackage{amsmath}
\usepackage{amsthm}
\usepackage{tikz}
\usepackage{url}
\usepackage{hyperref}
\usepackage{nicematrix}
\usepackage{breqn}
\newtheoremstyle{custom}
{6pt}% Space above
{3pt}% Space below
{\sl}% Body font
{}% Indent amount
{\bf}% Theorem head font
{.}% Punctuation after theorem head 
{.5em}% Space after theorem head
{}% Theorem head spec (can be left empty, meaning `normal')

\theoremstyle{custom}
\newtheorem{theorem}{Theorem}[section]

\newtheorem{conjecture}[theorem]{Conjecture}
\newtheorem{remark}[theorem]{Remark}

\newtheorem{statement}[theorem]{Statement}

\newtheorem{question}[theorem]{Question}

\title{On a polynomial positivity question of Collins, Dykema, and Torres-Ayala related to the Bessis-Moussa-Villani Conjecture}
\author{Edward D.~Kim, \sl{University of Wisconsin-La Crosse}\\
Joe Miller, \sl{University of Wisconsin-La Crosse}\\
Laura Zinnel, \sl{Iowa State University}}

\begin{document}

\maketitle

\begin{abstract}
We construct sum of squares certificates of non-negativity for two families of polynomials appearing as a variant by Collins, Dykema, and Torres-Ayala to H\"aegele's reformulation of a conjecture by Bessis, Moussa, and Villani.
\end{abstract}

\section{Introduction}

In 1975, Bessis, Moussa, and Villani (see~\cite{BessisMoussaVillani}) conjectured the following statement due to its potential use in quantum mechanical partition functions:
\begin{statement}\label{statement:bmv-original}
Given an $n \times n$ Hermitian matrix $A$ and an $n \times n$ positive semidefinite matrix $B$, the function $t \mapsto \operatorname{trace}(\exp(A-tB))$ is the Laplace transform of a positive measure supported in $[0, \infty)$.
\end{statement}
This was proved by Stahl (see~\cite{Stahl:ProofBMV}), published posthumously in 2013. A simplified proof (see~\cite{Eremenko}) is due to Eremenko. For a survey on the problem, see~\cite{Moussa}. In~\cite{LiebSeiringer}, Lieb and Seiringer gave the following reformulation of the statement:
\begin{conjecture}\label{conjecture:LiebSeiringer}
Given $n \times n$ positive semidefinite matrices $A$ and $B$, for all integers $m \geq r \geq 0$, the coefficient of $t^r$ in $\operatorname{trace}((A+tB)^m)$ is non-negative.
\end{conjecture}
Hillar and Johnson (see~\cite{HillarJohnson:n2m5}) proved that Conjecture~\ref{conjecture:LiebSeiringer} holds if $n=2$. In~\cite{HillarJohnson}, Hillar and Johnson proved Conjecture~\ref{conjecture:LiebSeiringer} also holds for the case $n=3$ and $m=6$.
\begin{theorem}[Hillar Descent Theorem, \cite{Hillar:Advances}]
Suppose that there exist integers $M$ and $R$ and $n \times n$ positive definite matrices $A$ and $B$ such that the coefficient of $t^R$ in $\operatorname{trace}((A+tB)^M)$ is negative. Then, for any $m \geq M$ and $r \geq R$ such that $m-r \geq M-R$, there exist $n \times n$ positive definite $A'$ and $B'$ such that the coefficient of $t^r$ in $\operatorname{trace}((A'+tB')^m)$ is negative.
\end{theorem}
Thus, it is enough to prove Conjecture~\ref{conjecture:LiebSeiringer} for infinitely many $m$. More specifically, it is enough to prove that the coefficient of $t^r$ in $\operatorname{trace}((A+tB)^m)$ is non-negative for an infinite sequence $(m_i,r_i)$ with: $m_i \rightarrow \infty$, $r_i \rightarrow \infty$, and $m_i-r_i \rightarrow \infty$.

H\"agele (see~\cite{Haegele}) proved Conjecture~\ref{conjecture:LiebSeiringer} holds for $m = 7$, and by the Hillar Descent Theorem it holds for $m \leq 7$. With a strengthened version of the Hillar Descent Theorem, Klep and Schweighofer (see~\cite{KS:SumsHermitianSquares}) proved Conjecture~\ref{conjecture:LiebSeiringer} holds for $m \leq 13$.

Positive results on Conjecture~\ref{conjecture:LiebSeiringer} are also known based on which coefficient of $t^r$ is examined. Burgdorf (see~\cite{Burgdorf}) proved Conjecture~\ref{conjecture:LiebSeiringer} holds for $r=4$. (That is, the fourth coefficient is always non-negative.) Due to the Hillar Descent Theorem, the conjecture holds for $r \leq 4$ and by symmetry for $r \geq m-4$. In~\cite{FleischhackFriedland}, Fleischhack and Friedland proved an asymptotic result: if $AB \not= 0$, then for fixed $A$, $B$, and $r$, there exists $m_0$ such that the coefficient of $t^r$ is non-negative for all $m \geq m_0$.

Let $S_{m,r}(A,B)$ be the sum of all words in the two non-commuting letters $A$ and $B$ with $r$ letters equal to $B$ and $m-r$ letters equal to $A$. We say that a polynomial $f$ in the ring $R$ of polynomials in non-commuting variables is cyclically equivalent to a sum of Hermitian squares if $f$ is a sum of commutators and Hermitian squares. It was hoped that Conjecture~\ref{conjecture:LiebSeiringer} could be proved by showing that $S_{m,r}(A,B)$ is a sum of commutators and Hermitian squares in the non-commuting variables $A$ and $B$: since $A$ and $B$ are positive semidefinite, there exist symmetric matrices $X$ and $Y$ such that $A=X^2$ and $B=Y^2$, and moreover, if $S_{m,r}(X^2,Y^2)$ is cyclically equivalent to a sum of Hermitian squares for given $(m,r)$, then Conjecture~\ref{conjecture:LiebSeiringer} is true for the case $(m,r)$. That is, $S_{m,r}(A,B)$ is cyclically equivalent to a sum of Hermitian squares implies the coefficient of $t^r$ in $\operatorname{trace}((A+tB)^m)$ is non-negative. % Observation 1.6 from CKP

In~\cite{Landweber}, Landweber and Speer prove $S_{m,r}(X^2,Y^2)$ is not cyclically equivalent to a sum of Hermitian squares when $m$ is odd and $5 \leq r \leq m-5$, or when $m$ is even and $r$ is odd and $3 \leq r \leq m-3$, or when $m \geq 13$ and $r=3$, or when $(m,r)=(9,3)$. When this is combined with the result of Collins, Dykema, and Torres-Ayala in~\cite{Collins} that $S_{m,r}(X^2,Y^2)$ is not cyclically equivalent to a sum of Hermitian squares when $m$ and $r$ are even and $6 \leq r \leq m-10$, this rules out the avenue of the approach of writing $S_{m,r}$ as a noncommutative sum of squares to prove Conjecture~\ref{conjecture:LiebSeiringer}. In~\cite{Cafuta} Cafuta, Klep, and Povh prove the then-remaining case, namely, $S_{16,8}(X^2,X^2)$ is not cyclically equivalent to a sum of Hermitian squares. In summary, $S_{m,r}(X^2,Y^2)$ is not cyclically equivalent to a sum of Hermitian squares when $m-6 \geq r \geq 6$ and $m \geq 16$.

Collins, Dykema, and Torres-Ayala (see~\cite{Collins}) modified the Lieb-Seiringer formulation given in Conjecture~\ref{conjecture:LiebSeiringer} in two ways, considering the following question:
\begin{conjecture}\label{conjecture:CDTA}
Given $n \times n$ symmetric matrices $A$ and $B$, for all $m \geq r$ where $m$ and $r$ are even, the coefficient of $t^r$ in $\operatorname{trace}((A+tB)^m)$ is non-negative.
\end{conjecture}
This conjecture is curiously neither stronger than nor weaker than Conjecture~\ref{conjecture:LiebSeiringer} as the statement relaxes the positive-semidefiniteness of $A$ and $B$ to symmetry, but strengthens the requirement on $m$ and $r$ by demanding evenness. However,~\cite{Collins} notes that a positive result to this conjecture for certain parameters $(m,r,n)$ implies that Conjecture~\ref{conjecture:LiebSeiringer} holds for the same parameters $(m,r,n)$, and thus together with Hillar's Descent Theorem, proving Conjecture~\ref{conjecture:CDTA} for an infinite sequence $(m_i,r_i)$ of special cases with $m_i \rightarrow \infty$, $r_i \rightarrow \infty$, and $m_i-r_i \rightarrow \infty$ would give an unexpected independent proof of Conjecture~\ref{conjecture:LiebSeiringer} and thus of Statement~\ref{statement:bmv-original}.

Previous work on Conjecture~\ref{conjecture:LiebSeiringer} involved determining whether $S_{m,r}(X^2,Y^2)$ is cyclically equivalent to a sum of Hermitian squares in two non-commutative variables, which was established positively for some $(m,r)$, but the results in~\cite{Collins, Landweber} gave a sufficient number of pairs $(m,r)$ for which this doesn't hold, precluding non-commutative sums of squares as an approach to prove Conjecture~\ref{conjecture:LiebSeiringer}. In Theorem~\ref{theorem:main-m4r2}, we give a proof -- independent of the proofs by Collins, Dykema, and Torres-Ayala in~\cite{Collins} or by Burgdorf et al.{} in~\cite{Burgdorf:TracialMoment, Burgdorf:CyclicEquivalenceChapter} -- of the special case of Conjecture~\ref{conjecture:CDTA} when $(m,r)=(4,2)$ by expressing the coefficient of $t^r$ in $\operatorname{trace}((A+tB)^m)$ as a sum of squares in the $2n + 2\binom{n}{2}$ commutative variables $a_{\{i,j\}}$ and $b_{\{i,j\}}$ which are the \emph{entries} of the $n \times n$ symmetric matrices $A$ and $B$.

Constructing a sum of squares certificate in $2n + 2\binom{n}{2}$ commutative variables admittedly reverts from the methodology of H\"agele in~\cite{Haegele} and subsequent articles. Moreover, due to the structure observed in proving Theorem~\ref{theorem:main-m4r2}, in Section~\ref{section:m4comments}, we present an even \emph{further} strengthening of Conjecture~\ref{conjecture:CDTA} and give evidence that the question merits further attention.

These remarks and methodology led us to the following: in Theorem~\ref{theorem:main-m8r4}, we give a proof of the special case of Conjecture~\ref{conjecture:CDTA} when $(m,r)=(8,4)$ for sufficiently small $n$.

Our results rely on a well-known (see, for instance~\cite{Laurent} and the references therein) equivalence, namely, a given multivariate polynomial is a sum of squares of polynomials (and hence non-negative) if and only if a certain semidefinite optimization problem is feasible.

\section{Background and Summary of Results}

Throughout, $[n]$ denotes $\{1,\dots,n\}$. Since $A$ and $B$ are symmetric matrices, we can denote the $(i,j)$-entry of $A$ and $B$ respectively by $d^0_{\{i,j\}}$ and $d^1_{\{i,j\}}$ respectively. Then, the $(i,j)$-entry of $(A+tB)^m$ is
\[\sum_{k_1 \in [n]} \sum_{k_2 \in [n]} \cdots \sum_{k_{m-1} \in [n]} \sum_{s_1+s_2+\dots+s_m=r} d^{s_1}_{\{i,k_1\}} d^{s_2}_{\{k_1,k_2\}} \cdots d^{s_m}_{\{k_{m-1},j\}} t^r,\]
where $s_1,\dots,s_m$ are in $\{0,1\}$.
Thus, for a fixed $r$, the coefficient of $t^r$ in the trace of $(A+tB)^m$ is 
\begin{equation}\label{equation:pre-necklace}
\sum_{k_0 \in [n]} \sum_{k_1 \in [n]} \sum_{k_2 \in [n]} \cdots \sum_{k_{m} \in [n]} \sum_{s_1+s_2+\dots+s_m=r} d^{s_1}_{\{k_0,k_1\}} d^{s_2}_{\{k_1,k_2\}} \cdots d^{s_m}_{\{k_{m-1},k_0\}}.
\end{equation}
The use of $s_1,\dots,s_m \in \{0,1\}$ together with the entries of $A$ and $B$ denoted by $d^0_{\{i,j\}}$ and $d^1_{\{i,j\}}$ is advantageous to describe expression~\eqref{equation:pre-necklace} unambiguously, but it will generally be more convenient for us to use $a_{\{i,j\}}$ and $b_{\{i,j\}}$ to denote the entries of $A$ and $B$.

The special case of Conjecture~\ref{conjecture:CDTA} where $r=0$ is trivial as the coefficient of $t^0$ in $\operatorname{trace}((A+tB)^m)$ is
\[\sum_{i,j \in [n]}\left(\sum_{k_1 \in [n]} \sum_{k_2 \in [n]} \cdots \sum_{k_{m/2-1} \in [n]} a_{\{i,k_1\}} a_{\{k_1,k_2\}} \cdots a_{\{k_{m/2-1},j\}}\right)^2\]
and by symmetry, the case when $r=m$ simply uses the $b$ variables instead of the $a$ variables. We now examine the first non-trivial case, where $(m,r)=(4,2)$. We prove the following theorem in Section~\ref{section:m4r2}:
\begin{theorem}\label{theorem:main-m4r2}
Given $n \times n$ symmetric matrices $A$ and $B$, the coefficient of $t^2$ in $\operatorname{trace}((A+tB)^4)$ is non-negative. Moreover, the coefficient of $t^2$ in $\operatorname{trace}((A+tB)^4)$ is a sum of squares of polynomials in the variables $a_{\{i,j\}}$ and $b_{\{i,j\}}$: there exist positive semidefinite matrices $Q_1$ and $Q_2$, and vectors $z_1$ and $z_{2,(i,j)}$ such that the coefficient of $t^2$ in $\operatorname{trace}((A+tB)^4)$ is equal to 
\begin{equation}\label{equation:big-equation}
z_1^TQ_1z_1 + \sum_{1 \leq i < j \leq n} z_{2,(i,j)}^TQ_2z_{2,(i,j)}.
\end{equation}
\end{theorem}

\begin{figure}[hbt]
\begin{center}
\begin{tikzpicture}
\draw(.5,-.5) -- (1.5,-.5) -- (1.5,.5) -- (-1.5,.5) -- (-1.5,-.5) -- (-.5,-.5); %\draw(1.5,.5) -- (-1.5,.5) -- (-1.5,-.5) -- (1.5,-.5) -- cycle;
\node at (0,-.5) {$\dots$};
\draw[fill=white] (1.5,.5) circle (.2);
\draw[fill=white] (.5,.5) circle (.2);
\draw[fill=white] (-.5,.5) circle (.2);
\draw[fill=white] (-1.5,.5) circle (.2);
\draw[fill=white] (-1.5,-.5) circle (.2);
\draw[fill=white] (-.5,-.5) circle (.2);
\draw[fill=white] (.5,-.5) circle (.2);
\draw[fill=white] (1.5,-.5) circle (.2);
\node at (1.5,.5) {$s_1$}; 
\node at (.5,.5) {$s_2$}; 
\node at (-.5,.5) {$s_3$}; 
\node at (-1.5,.5) {$s_4$}; 
\node at (-1.5,-.5) {$s_5$}; 
\node at (-.5,-.5) {$s_6$}; 
\node at (.5,-.5) {{\tiny$s_{m-1}$}}; 
\node at (1.5,-.5) {$s_m$}; 
\node at (1,.7) {$k_1$}; 
\node at (0,.7) {$k_2$}; 
\node at (-1,.7) {$k_3$}; 
\node at (-1.7,0) {$k_4$}; 
\node at (-1,-.7) {$k_5$}; 
\node at (1,-.7) {{\tiny$k_{m-1}$}}; 
\node at (1.7,0) {$k_0$}; 
\end{tikzpicture}
\caption{An $(m,r,n)$-necklace has exactly $r$ of $s_1,\dots,s_m$ as $1$'s and $m-r$ as $0$'s}\label{figure:generic-mrn-necklace}
\end{center}
\end{figure}
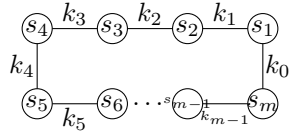

Motivated by expression~\eqref{equation:pre-necklace}, we define an $(m,r,n)$-necklace (or simply a necklace if the parameters are clear from context) to be a vertex-labeled edge-labeled $m$-cycle (see Figure~\ref{figure:generic-mrn-necklace}) where each vertex label is in $[n]$, each edge label is in $\{0,1\}$, and precisely $r$ of the edge labels are $1$s. When~\eqref{equation:pre-necklace} is expanded without any grouping of terms, the set of terms is in one-to-one correspondence with the set of $(m,r,n)$-necklaces. It will be convenient to use $a$'s and $b$'s instead of $0$'s and $1$'s for the vertex labels $s_1,\dots,s_m$.

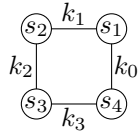
\begin{figure}[hbt]
\begin{center}
\begin{tikzpicture}
\draw(.5,.5) -- (-.5,.5) -- (-.5,-.5) -- (.5,-.5) -- cycle;
\draw[fill=white] (.5,.5) circle (.2);
\draw[fill=white] (-.5,.5) circle (.2);
\draw[fill=white] (-.5,-.5) circle (.2);
\draw[fill=white] (.5,-.5) circle (.2);
\node at (.5,.5) {$s_1$}; % northeast
\node at (-.5,.5) {$s_2$}; % northwest
\node at (-.5,-.5) {$s_3$}; % southwest
\node at (.5,-.5) {$s_4$}; % southeast
\node at (0,.7) {$k_1$}; % north
\node at (-.7,0) {$k_2$}; % west
\node at (0,-.7) {$k_3$}; % south
\node at (.7,0) {$k_0$}; % east
\end{tikzpicture}
\caption{In a $(4,2,n)$-necklace, two of $s_1,\dots,s_4$ are $a$'s and two are $b$'s}\label{figure:generic-m4r2}
\end{center}
\end{figure}

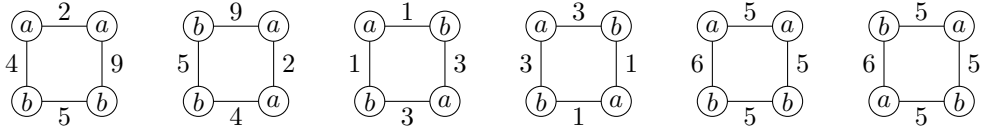
\begin{figure}[hbt]
\begin{center}
\begin{tikzpicture}
\draw(.5,.5) -- (-.5,.5) -- (-.5,-.5) -- (.5,-.5) -- cycle;
\draw[fill=white] (.5,.5) circle (.2);
\draw[fill=white] (-.5,.5) circle (.2);
\draw[fill=white] (-.5,-.5) circle (.2);
\draw[fill=white] (.5,-.5) circle (.2);
\node at (.5,.5) {$a$}; % northeast
\node at (-.5,.5) {$a$}; % northwest
\node at (-.5,-.5) {$b$}; % southwest
\node at (.5,-.5) {$b$}; % southeast
\node at (0,.7) {$2$}; % north
\node at (-.7,0) {$4$}; % west
\node at (0,-.7) {$5$}; % south
\node at (.7,0) {$9$}; % east
\end{tikzpicture}
\quad
\begin{tikzpicture}
\draw(.5,.5) -- (-.5,.5) -- (-.5,-.5) -- (.5,-.5) -- cycle;
\draw[fill=white] (.5,.5) circle (.2);
\draw[fill=white] (-.5,.5) circle (.2);
\draw[fill=white] (-.5,-.5) circle (.2);
\draw[fill=white] (.5,-.5) circle (.2);
\node at (.5,.5) {$a$}; % northeast
\node at (-.5,.5) {$b$}; % northwest
\node at (-.5,-.5) {$b$}; % southwest
\node at (.5,-.5) {$a$}; % southeast
\node at (0,.7) {$9$}; % north
\node at (-.7,0) {$5$}; % west
\node at (0,-.7) {$4$}; % south
\node at (.7,0) {$2$}; % east
\end{tikzpicture}
\quad
\begin{tikzpicture}
\draw(.5,.5) -- (-.5,.5) -- (-.5,-.5) -- (.5,-.5) -- cycle;
\draw[fill=white] (.5,.5) circle (.2);
\draw[fill=white] (-.5,.5) circle (.2);
\draw[fill=white] (-.5,-.5) circle (.2);
\draw[fill=white] (.5,-.5) circle (.2);
\node at (.5,.5) {$b$}; % northeast
\node at (-.5,.5) {$a$}; % northwest
\node at (-.5,-.5) {$b$}; % southwest
\node at (.5,-.5) {$a$}; % southeast
\node at (0,.7) {$1$}; % north
\node at (-.7,0) {$1$}; % west
\node at (0,-.7) {$3$}; % south
\node at (.7,0) {$3$}; % east
\end{tikzpicture}
\quad
\begin{tikzpicture}
\draw(.5,.5) -- (-.5,.5) -- (-.5,-.5) -- (.5,-.5) -- cycle;
\draw[fill=white] (.5,.5) circle (.2);
\draw[fill=white] (-.5,.5) circle (.2);
\draw[fill=white] (-.5,-.5) circle (.2);
\draw[fill=white] (.5,-.5) circle (.2);
\node at (.5,.5) {$b$}; % northeast
\node at (-.5,.5) {$a$}; % northwest
\node at (-.5,-.5) {$b$}; % southwest
\node at (.5,-.5) {$a$}; % southeast
\node at (0,.7) {$3$}; % north
\node at (-.7,0) {$3$}; % west
\node at (0,-.7) {$1$}; % south
\node at (.7,0) {$1$}; % east
\end{tikzpicture}
\quad
\begin{tikzpicture}
\draw(.5,.5) -- (-.5,.5) -- (-.5,-.5) -- (.5,-.5) -- cycle;
\draw[fill=white] (.5,.5) circle (.2);
\draw[fill=white] (-.5,.5) circle (.2);
\draw[fill=white] (-.5,-.5) circle (.2);
\draw[fill=white] (.5,-.5) circle (.2);
\node at (.5,.5) {$a$}; % northeast
\node at (-.5,.5) {$a$}; % northwest
\node at (-.5,-.5) {$b$}; % southwest
\node at (.5,-.5) {$b$}; % southeast
\node at (0,.7) {$5$}; % north
\node at (-.7,0) {$6$}; % west
\node at (0,-.7) {$5$}; % south
\node at (.7,0) {$5$}; % east
\end{tikzpicture}
\quad
\begin{tikzpicture}
\draw(.5,.5) -- (-.5,.5) -- (-.5,-.5) -- (.5,-.5) -- cycle;
\draw[fill=white] (.5,.5) circle (.2);
\draw[fill=white] (-.5,.5) circle (.2);
\draw[fill=white] (-.5,-.5) circle (.2);
\draw[fill=white] (.5,-.5) circle (.2);
\node at (.5,.5) {$a$}; % northeast
\node at (-.5,.5) {$b$}; % northwest
\node at (-.5,-.5) {$a$}; % southwest
\node at (.5,-.5) {$b$}; % southeast
\node at (0,.7) {$5$}; % north
\node at (-.7,0) {$6$}; % west
\node at (0,-.7) {$5$}; % south
\node at (.7,0) {$5$}; % east
\end{tikzpicture}
\caption{Some $(4,2,n)$-necklaces}\label{figure:specific-m4r2}
\end{center}
\end{figure}

To prove Theorem~\ref{theorem:main-m4r2}, every $(4,2,n)$-necklace (whose general form is depicted in Figure~\ref{figure:generic-m4r2}, with precisely two of $s_1,\dots,s_4$ being $a$'s and two being $b$'s, and $k_0,\dots,k_3 \in [n]$) will be considered exactly once in our construction of a sum of squares expression.
Each necklace corresponds to a monomial in~\eqref{equation:pre-necklace}. For example, the first $(4,2,n)$-necklace in Figure~\ref{figure:specific-m4r2} relates to the monomial $a_{\{9,2\}}a_{\{2,4\}}b_{\{4,5\}}b_{\{5,9\}}$, the second relates to the monomial $a_{\{2,9\}}b_{\{9,5\}}b_{\{5,4\}}a_{\{4,2\}}$, the third relates to the monomial $b_{\{3,1\}}a_{\{1,1\}}b_{\{1,3\}}a_{\{3,3\}}$, the fourth necklace relates to the monomial $b_{\{1,3\}}a_{\{3,3\}}b_{\{3,1\}}a_{\{1,1\}}$, the fifth necklace relates to the monomial $a_{\{5,5\}}a_{\{5,6\}}b_{\{6,5\}}b_{\{5,5\}}$, and the sixth necklace relates to the monomial $a_{\{5,5\}}b_{\{5,6\}}a_{\{6,5\}}b_{\{5,5\}}$.

Since $A$ and $B$ are symmetric, we note that the first and second necklaces in Figure~\ref{figure:specific-m4r2} relate to the same monomial, the third and fourth necklaces relate to the same monomial, and the fifth and sixth necklaces relate to the same monomial. While the first and second necklaces in Figure~\ref{figure:specific-m4r2} can be obtained from each other via dihedral group action, and the same is true for obtaining the third and fourth necklaces from each other, we note that the fifth and sixth necklaces relate to the same monomial but cannot be obtained from each other using an action from the dihedral group of a square.

We prove the following theorem in Section~\ref{section:m8r4} which provides a partial resolution to a question of Collins, Dykema, and Torres-Ayala:
\begin{theorem}\label{theorem:main-m8r4}
Given $n \times n$ symmetric matrices $A$ and $B$ with $n \leq 5$, the coefficient of $t^4$ in $\operatorname{trace}((A+tB)^8)$ is non-negative. Moreover, the coefficient of $t^4$ in $\operatorname{trace}((A+tB)^8)$ is a sum of squares of polynomials in the variables $a_{\{i,j\}}$ and $b_{\{i,j\}}$: there exist positive semidefinite matrices $Q_1$ and $Q_2$ and $Q_3$, and vectors $z_1$ and $z_2$ and $z_{3,(i,j)}$ such that the coefficient of $t^4$ in $\operatorname{trace}((A+tB)^8)$ is equal to 
\begin{equation}\label{equation:big-equation-m8}
z_1^TQ_1z_1 + z_2^TQ_2z_2 + \sum_{1 \leq i < j \leq n} z_{3,(i,j)}^TQ_3z_{3,(i,j)}.
\end{equation}
\end{theorem}

\section{Proof of Theorem~\ref{theorem:main-m4r2}}\label{section:m4r2}

Theorem~\ref{theorem:main-m4r2} follows by the following lemma and proving the matrices $Q_1$ and $Q_2$ appearing are positive semidefinite. Since the coefficient of $t^2$ in $\operatorname{trace}((A+tB)^4)$ is equivalent to enumerating $(4,2,n)$-necklaces, the proof  will consider all $(4,2,n)$-necklaces.

In Section~\ref{section:m4defineQz}, we define the relevant SDP matrices and monomial vectors, and in Section~\ref{section:m4exampleQz} give an example. In Section~\ref{section:m4provePSD}, we prove positive semidefiniteness, and in Section~\ref{section:m4-equal-expressions} we prove that our sum of squares representation is indeed for the coefficient of $t^2$ in $\operatorname{trace}((A+tB)^4)$.

\subsection{Definitions of relevant matrices and monomial vectors}\label{section:m4defineQz}

Let $A$ and $B$ be $n \times n$ symmetric matrices whose entries are denoted $a_{\{i,j\}}$ and $b_{\{i,j\}}$, respectively.

Let $z_1$ be the vector of size $n+\binom{n}{2}$ indexed by all subsets of $[n]$ with cardinality $1$ or $2$, whose $\{i,j\}$-entry is $a_{\{i,j\}}b_{\{i,j\}}$. (Note that $i=j$ in the case where the indexing set is a singleton.)

Let $Q_1$ be the $(n+\binom{n}{2}) \times (n+\binom{n}{2})$ matrix whose rows and columns are indexed by subsets of $[n]$ with cardinality $1$ or $2$, whose $(\{i,j\}, \{i',j'\})$-entry is $6 | \{i,j\} \cap \{i',j'\} |$. %$(\{i,j\},\{k,\ell\})$-entry is the scaled cardinality $6 \big| \{i,j\} \cap \{k,\ell\}\big|$.

For $1 \leq i < j \leq n$, let $z_{2,(i,j)}$ be the vector
\[ 
\left(
\begin{array}{c}
a_{\{i,1\}} b_{\{j,1\}} \\
a_{\{i,2\}} b_{\{j,2\}} \\
\vdots\\
a_{\{i,n\}} b_{\{j,n\}} \\ \hline
a_{\{j,1\}} b_{\{i,1\}} \\
a_{\{j,2\}} b_{\{i,2\}} \\
\vdots\\
a_{\{j,n\}} b_{\{i,n\}}
\end{array}
\right)\]
of size $2n$.

Let $Q_2$ be the $(2n) \times (2n)$ matrix
\[\left( \begin{array}{c|c} 4 J_n & 2 J_n \\\hline 2 J_n & 4 J_n \end{array}\right),\]
where $J_n$ denotes the $n \times n$ all $1$'s matrix.

\begin{remark}
The expression given in~\eqref{equation:big-equation} is admittedly asymmetric due to $i < j$ appearing in the sum. The expression can be symmetrized to 
\[z_1^TQ_1z_1 + \frac12\sum_{i\not=j} z_{2,(i,j)}^TQ_2z_{2,(i,j)}\]
but this would unnecessarily complicate our exposition in the proof accounting for necklaces and matrix entries. In Section~\ref{section:m4-equal-expressions}, $i$ and $j$ will be used as edge labels on necklaces and it will not necessarily be the case that $i < j$ in the necklace diagrams.
\end{remark}

\subsection{Examples of relevant matrices and monomial vectors}\label{section:m4exampleQz}

To make the descriptions of the matrices and vectors concrete, we present these matrices and vectors when $n=3$. Then, the coefficient of $t^2$ in $\operatorname{trace}((A+tB)^4)$ is
\[z_1^T Q_1 z_1 
+ z_{2,(1,2)}^T Q_2 z_{2,(1,2)}
+ z_{2,(1,3)}^T Q_2 z_{2,(1,3)}
+ z_{2,(2,3)}^T Q_2 z_{2,(2,3)}
,\] where $Q_1$ is
\begin{equation*}
%\left(\begin{array}{ccc|ccc}
%6&0&0&6&6&0\\
%0&6&0&6&0&6\\
%0&0&6&0&6&6\\\hline
%6&6&0&12&6&6\\
%6&0&6&6&12&6\\
%0&6&6&6&6&12
%\end{array}\right)
\begin{pNiceArray}{ccc|ccc}[first-row,first-col]
& \{1\} & \{2\} & \{3\} & \{1,2\} & \{1,3\} & \{2,3\} \\
\{1\} &6&0&0&6&6&0\\
\{2\} &0&6&0&6&0&6\\
\{3\} &0&0&6&0&6&6\\\hline
\{1,2\} &6&6&0&12&6&6\\
\{1,3\} &6&0&6&6&12&6\\
\{2,3\} &0&6&6&6&6&12
\end{pNiceArray}
\end{equation*}
and $z_1$ is
\begin{equation*}
%z_1 = 
%\left(\begin{array}{c}
%a_{\{1,1\}} b_{\{1,1\}} \\
%a_{\{2,2\}} b_{\{2,2\}} \\
%a_{\{3,3\}} b_{\{3,3\}} \\ \hline
%a_{\{1,2\}} b_{\{1,2\}} \\
%a_{\{1,3\}} b_{\{1,3\}} \\
%a_{\{2,3\}} b_{\{2,3\}}
%\end{array}\right)
\begin{pNiceArray}{c}[first-col]
\{1\} & a_{\{1,1\}} b_{\{1,1\}} \\
\{2\} & a_{\{2,2\}} b_{\{2,2\}} \\
\{3\} & a_{\{3,3\}} b_{\{3,3\}} \\ \hline
\{1,2\} & a_{\{1,2\}} b_{\{1,2\}} \\
\{1,3\} & a_{\{1,3\}} b_{\{1,3\}} \\
\{2,3\} & a_{\{2,3\}} b_{\{2,3\}}
\end{pNiceArray}
\end{equation*}
and
\begin{equation*}
Q_2 =\left(
\begin{array}{ccc|ccc}
4&4&4&2&2&2\\
4&4&4&2&2&2\\
4&4&4&2&2&2\\\hline
2&2&2&4&4&4\\
2&2&2&4&4&4\\
2&2&2&4&4&4
\end{array}\right),
\end{equation*}
and
\begin{equation*}
z_{2,(1,2)} = \left(\begin{array}{c}
a_{\{1,1\}} b_{\{2,1\}}\\
a_{\{1,2\}} b_{\{2,2\}}\\
a_{\{1,3\}} b_{\{2,3\}}\\ \hline
a_{\{2,1\}} b_{\{1,1\}}\\
a_{\{2,2\}} b_{\{1,2\}}\\
a_{\{2,3\}} b_{\{1,3\}}
\end{array}\right),
\quad
z_{2,(1,3)} = \left(\begin{array}{c}
a_{\{1,1\}} b_{\{3,1\}}\\
a_{\{1,2\}} b_{\{3,2\}}\\
a_{\{1,3\}} b_{\{3,3\}}\\ \hline
a_{\{3,1\}} b_{\{1,1\}}\\
a_{\{3,2\}} b_{\{1,2\}}\\
a_{\{3,3\}} b_{\{1,3\}}
\end{array}\right),
\quad
z_{2,(2,3)} = \left(\begin{array}{c}
a_{\{2,1\}} b_{\{3,1\}}\\
a_{\{2,2\}} b_{\{3,2\}}\\
a_{\{2,3\}} b_{\{3,3\}}\\ \hline
a_{\{3,1\}} b_{\{2,1\}}\\
a_{\{3,2\}} b_{\{2,2\}}\\
a_{\{3,3\}} b_{\{2,3\}}
\end{array}\right),
\end{equation*}

\subsection{Proofs of positive semidefiniteness}\label{section:m4provePSD}

Recall that $Q_1$ is the matrix whose rows and columns are indexed by subsets of $[n]$ of cardinality $1$ or $2$, and whose $(\{i,j\}, \{i',j'\})$-entry is $6|\{i,j\} \cap \{i',j'\}|$. Let $U$ be the $n \times (n+ \binom{n}{2})$ matrix with rows are indexed by $[n]$ and whose columns are indexed by subsets of $[n]$ of cardinality $1$ or $2$, and the $(k,\{i,j\})$-entry of $U$ is $1$ if $k \in \{i,j\}$ and $0$ otherwise. For example, when $n=3$, then $U$ is
\[
\begin{pNiceArray}{ccc|ccc}[first-row,first-col]
& \{1\} & \{2\} & \{3\} & \{1,2\} & \{1,3\} & \{2,3\} \\
1 &1&0&0&1&1&0\\
2 &0&1&0&1&0&1\\
3 &0&0&1&0&1&1
\end{pNiceArray}.
\]
Then $Q_1 = 6\,U^TU$, proving $Q_1$ is positive semidefinite. Note $Q_1$ is a scaled example of a matrix recording the cardinalities of set intersections as described in Section~1.3 of~\cite{BrualdiRyser}. Our matrices are reminiscent of, but not identical to, the combinatorial matrices in Knuth's note~\cite{Knuth:cm}.

Note that $Q_2$ can be easily shown to be a sum of positive semidefinite matrices. For an alternate but succinct proof, $Q_2$ is the tensor product of the matrix
\[\left(\begin{array}{cc}4&2\\2&4\end{array}\right)\]
and the $n \times n$ all $1$'s matrix. It is well-known (see for instance Theorem 9.1.5 in~\cite{BrualdiCvetkovic}) that the tensor product preserves positive semidefiniteness.

\subsection{Notes prior to necklace counting}

We have now shown that~\eqref{equation:big-equation} is a sum of squares expression. Prior to proving this expression is indeed equal to the coefficient of $t^2$ in $\operatorname{trace}((A+tB)^4)$, we introduce terminology to standardize language we will use to describe the entries of $Q_1$ and $Q_2$.

The matrix $Q_1$ contains an upper left block of $6I$ where $I$ is the $n\times n$ identity matrix. We will call this the upper left block of $Q_1$. The $\binom{n}{2} \times n$ block below the upper left block of $Q_1$ and the $n\times \binom{n}{2}$ right of the upper left block of $Q_1$ will be called the lower left and upper right blocks of $Q_1$, respectively. Lastly, the square $\binom{n}{2} \times \binom{n}{2}$ block on the bottom right will be called lower right block of $Q_1$.

The matrix $Q_2$ consists of upper-left and lower-right $n \times n$ blocks of all $4$'s. We refer to these as the {\bf diagonal blocks}. The upper-right and lower-left $n \times n$ blocks of all $2$'s are referred to as the {\bf non-diagonal blocks}. Since $Q_2$ is $(2n) \times (2n)$, it is technically accurate to refer to the $(n+u,v)$-entry of $Q_2$, but it will be more convenient to refer to this non-diagonal block entry as the $(u,v)$-entry of the lower left block of $Q_2$. Similarly, the $(n+u,n+v)$-entry of $Q_2$ is also called the $(u,v)$-entry of the lower right diagonal block of $Q_2$.

Based on our description of the blocks of $Q_2$, it is natural to describe the vector $z_{2,(i,j)}$ of size $2n$ as having an upper half block and a lower half block. The $k$th entry of the upper half block of $z_{2,(i,j)}$ is $a_{\{i,k\}}b_{\{j,k\}}$ while the $k$th entry of the lower half block is $a_{\{j,k\}}b_{\{i,k\}}$.

\begin{figure}[hbt]
\begin{center}
\begin{tikzpicture}
\draw(.5,.5) -- (-.5,.5) -- (-.5,-.5) -- (.5,-.5) -- cycle;
\draw[fill=white] (.5,.5) circle (.2);
\draw[fill=white] (-.5,.5) circle (.2);
\draw[fill=white] (-.5,-.5) circle (.2);
\draw[fill=white] (.5,-.5) circle (.2);
\node at (.5,.5) {$a$}; % northeast
\node at (-.5,.5) {$b$}; % northwest
\node at (-.5,-.5) {$a$}; % southwest
\node at (.5,-.5) {$b$}; % southeast
\node at (0,.7) {$1$}; % north
\node at (-.7,0) {$2$}; % west
\node at (0,-.7) {$3$}; % south
\node at (.7,0) {$5$}; % east
\end{tikzpicture}
\quad
\begin{tikzpicture}
\draw(.5,.5) -- (-.5,.5) -- (-.5,-.5) -- (.5,-.5) -- cycle;
\draw[fill=white] (.5,.5) circle (.2);
\draw[fill=white] (-.5,.5) circle (.2);
\draw[fill=white] (-.5,-.5) circle (.2);
\draw[fill=white] (.5,-.5) circle (.2);
\node at (.5,.5) {$b$}; % northeast
\node at (-.5,.5) {$a$}; % northwest
\node at (-.5,-.5) {$b$}; % southwest
\node at (.5,-.5) {$a$}; % southeast
\node at (0,.7) {$1$}; % north
\node at (-.7,0) {$5$}; % west
\node at (0,-.7) {$3$}; % south
\node at (.7,0) {$2$}; % east
\end{tikzpicture}
\caption{Two $(4,2,n)$-necklaces counted by an entry in an instance of $Q_2$}\label{figure:for-accounting-language}
\end{center}
\end{figure}
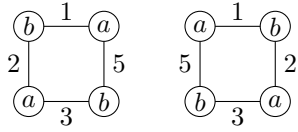

In Section~\ref{section:m4-equal-expressions} we account for all necklaces by specifying whether a given necklace is counted in $z_1^TQ_1z_1$ or $\sum z_{2,(i,j)}^TQ_2z_{2,(i,j)}$. More specifically, we will declare whether a given necklace is counted in $Q_1$ or counted in a copy of $Q_2$, where each entry of $Q_1$ and each entry of each copy of $Q_2$ is mentioned exactly once. Furthermore, due to the uniformity in $Q_2$, we will state whether a given necklaces is counted in a diagonal or non-diagonal block of a $Q_2$ matrix. To facilitate this process of considering each $(4,2,n)$-necklace exactly once and considering each entry of the $Q$ matrices exactly once, it will be helpful to have standard language in our accounting process in the next section: for example, we say the two necklaces depicted in Figure~\ref{figure:for-accounting-language} are counted by the $2$ appearing as the $(5,2)$-entry in the upper-right block of the $(1,3)$ copy of $Q_2$. 
Using the entries of the matrices appearing in~\eqref{equation:big-equation} to describe the number of necklaces satisfying specified conditions (which in turn counts the multiplicities of the represented monomials in~\eqref{equation:pre-necklace}) is justified since $z^TQz$ can also be written
\[\sum_{x,y} q_{x,y}z_xz_y.\] 
Moreover, counting in this manner is only possible since $Q_1$ and $Q_2$ are integer-valued matrices.

\subsection{Proof of relevance to~\eqref{equation:big-equation}}\label{section:m4-equal-expressions}

In this section, we prove the coefficient of $t^2$ in $\operatorname{trace}((A+tB)^4)$ is equal to the expression
\[z_1^TQ_1z_1 + \sum_{1 \leq i < j \leq n} z_{2,(i,j)}^TQ_2z_{2,(i,j)}\] 
from~\eqref{equation:big-equation} by considering the set of all $(4,2,r)$-necklaces and classifying them according to their structure. An edge in a necklace is called \emph{balanced} if its endpoints are a vertex labeled $a$ and a vertex labeled $b$, and is called \emph{unbalanced} otherwise. 

For example, in the first necklace in Figure~\ref{figure:specific-m4r2}, the left and right edges (with $4$ and $9$ labels) are balanced edges while the top and bottom edges are unbalanced. In the last necklace of Figure~\ref{figure:specific-m4r2}, the vertex labels alternate $a$'s and $b$'s, thus all four edges are balanced. 

Using the notation established in Figure~\ref{figure:generic-m4r2}, we partition the set of all $(4,2,n)$-necklaces into four mutually-distinct collections:
\begin{enumerate}
\item The $a$'s are consecutive and $b$'s are consecutive, and the two unbalanced edges have mismatched labels while the two balanced edges have matching labels
\item Otherwise, $k_1=k_3$, or $k_2=k_0$ and $s_2 \not= s_4$
\item Otherwise, $s_2 \not= s_4$
\item Otherwise.
\end{enumerate}
The four collections are clearly mutually-distinct since the definitions of subsequent collections explicitly mention ``otherwise'' to the point where one could verify this classification scheme computationally using nested if-else structures. These four collections can equivalently be described without subsequent collections relying on prior collections as follows:
\begin{enumerate}
\item The $a$'s are consecutive and $b$'s are consecutive, and the two unbalanced edges have mismatched labels while the two balanced edges have matching labels. (See Figure~\ref{figure:q2-diagonal}.) The four such necklaces satisfying these conditions (all obtained as rotations of each other) are accounted for by the $4$ appearing as the $(k,k)$-entry in the upper-left or lower-right block of a copy of $Q_2$. More specifically, if $i<j$, then these necklaces are counted by the $4$ appearing as the $(k,k)$-entry of the upper-left block of the $(i,j)$ copy of $Q_2$, but if $j<i$, then these necklaces are counted by the $4$ appearing as the $(k,k)$-entry of the lower-right block of the $(j,i)$ copy of $Q_2$.

\begin{figure}[hbt]
\begin{center}
\begin{tikzpicture}
\draw(.5,.5) -- (-.5,.5) -- (-.5,-.5) -- (.5,-.5) -- cycle;
\draw[fill=white] (.5,.5) circle (.2);
\draw[fill=white] (-.5,.5) circle (.2);
\draw[fill=white] (-.5,-.5) circle (.2);
\draw[fill=white] (.5,-.5) circle (.2);
\node at (.5,.5) {$a$}; % northeast
\node at (-.5,.5) {$a$}; % northwest
\node at (-.5,-.5) {$b$}; % southwest
\node at (.5,-.5) {$b$}; % southeast
\node at (0,.7) {$i$}; % north
\node at (-.7,0) {$k$}; % west
\node at (0,-.7) {$j$}; % south
\node at (.7,0) {$k$}; % east
\end{tikzpicture}
\begin{tikzpicture}
\draw(.5,.5) -- (-.5,.5) -- (-.5,-.5) -- (.5,-.5) -- cycle;
\draw[fill=white] (.5,.5) circle (.2);
\draw[fill=white] (-.5,.5) circle (.2);
\draw[fill=white] (-.5,-.5) circle (.2);
\draw[fill=white] (.5,-.5) circle (.2);
\node at (.5,.5) {$b$}; % northeast
\node at (-.5,.5) {$a$}; % northwest
\node at (-.5,-.5) {$a$}; % southwest
\node at (.5,-.5) {$b$}; % southeast
\node at (0,.7) {$k$}; % north
\node at (-.7,0) {$i$}; % west
\node at (0,-.7) {$k$}; % south
\node at (.7,0) {$j$}; % east
\end{tikzpicture}
\begin{tikzpicture}
\draw(.5,.5) -- (-.5,.5) -- (-.5,-.5) -- (.5,-.5) -- cycle;
\draw[fill=white] (.5,.5) circle (.2);
\draw[fill=white] (-.5,.5) circle (.2);
\draw[fill=white] (-.5,-.5) circle (.2);
\draw[fill=white] (.5,-.5) circle (.2);
\node at (.5,.5) {$b$}; % northeast
\node at (-.5,.5) {$b$}; % northwest
\node at (-.5,-.5) {$a$}; % southwest
\node at (.5,-.5) {$a$}; % southeast
\node at (0,.7) {$j$}; % north
\node at (-.7,0) {$k$}; % west
\node at (0,-.7) {$i$}; % south
\node at (.7,0) {$k$}; % east
\end{tikzpicture}
\begin{tikzpicture}
\draw(.5,.5) -- (-.5,.5) -- (-.5,-.5) -- (.5,-.5) -- cycle;
\draw[fill=white] (.5,.5) circle (.2);
\draw[fill=white] (-.5,.5) circle (.2);
\draw[fill=white] (-.5,-.5) circle (.2);
\draw[fill=white] (.5,-.5) circle (.2);
\node at (.5,.5) {$a$}; % northeast
\node at (-.5,.5) {$b$}; % northwest
\node at (-.5,-.5) {$b$}; % southwest
\node at (.5,-.5) {$a$}; % southeast
\node at (0,.7) {$k$}; % north
\node at (-.7,0) {$j$}; % west
\node at (0,-.7) {$k$}; % south
\node at (.7,0) {$i$}; % east
\end{tikzpicture}
\caption{necklaces with consecutive vertex labels satisfying $i \not= j$}\label{figure:q2-diagonal}
\end{center}
\end{figure}
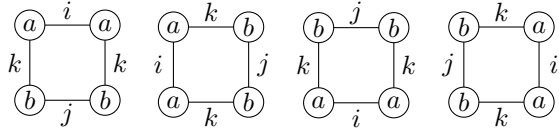

If $i < j$, the $k$th entry in the upper half of the vector $z_{2,(i,j)}$ is $a_{\{i,k\}}b_{\{j,k\}}$, thus the $(k,k)$-entry of the upper-left block is multiplied by $(a_{\{i,k\}}b_{\{j,k\}})^2$ in~\eqref{equation:big-equation}, a monomial represented by the necklaces in Figure~\ref{figure:q2-diagonal}. If $j < i$ then the $k$th entry in the lower half of the vector $z_{2,(j,i)}$ is $a_{\{i,k\}}b_{\{j,k\}}$, thus the $(k,k)$-entry of the lower-right block is multiplied by $(a_{\{i,k\}}b_{\{j,k\}})^2$ in~\eqref{equation:big-equation}, a monomial represented by the necklaces in Figure~\ref{figure:q2-diagonal}.

We now describe the remaining necklaces. That is, to avoid double counting, we summarize the set of necklaces not belonging to the first collection. Consider the necklaces which are not depicted in Figure~\ref{figure:q2-diagonal}. Among these necklaces, if the $a$'s and $b$'s do not alternate, the necklace must be one of the eight types of necklaces shown in Figure~\ref{figure:q2-diagonal-negation}, where unmarked edges may take on any label in $[n]$. The first row of Figure~\ref{figure:q2-diagonal-negation} considers the corresponding configuration of $a$'s and $b$'s from Figure~\ref{figure:q2-diagonal} where the $i \not= j$ condition fails, while the second row of Figure~\ref{figure:q2-diagonal-negation} considers the corresponding configuration of $a$'s and $b$'s from Figure~\ref{figure:q2-diagonal} depicts the condition of the matching edge labels in Figure~\ref{figure:q2-diagonal} now instead being mismatched.
\begin{figure}[hbt]
\begin{center}
\begin{tikzpicture}
\draw(.5,.5) -- (-.5,.5) -- (-.5,-.5) -- (.5,-.5) -- cycle;
\draw[fill=white] (.5,.5) circle (.2);
\draw[fill=white] (-.5,.5) circle (.2);
\draw[fill=white] (-.5,-.5) circle (.2);
\draw[fill=white] (.5,-.5) circle (.2);
\node at (.5,.5) {$a$}; % northeast
\node at (-.5,.5) {$a$}; % northwest
\node at (-.5,-.5) {$b$}; % southwest
\node at (.5,-.5) {$b$}; % southeast
\node at (0,.7) {$k$}; % north
\node at (-.7,0) {}; % west
\node at (0,-.7) {$k$}; % south
\node at (.7,0) {}; % east
\end{tikzpicture}
\begin{tikzpicture}
\draw(.5,.5) -- (-.5,.5) -- (-.5,-.5) -- (.5,-.5) -- cycle;
\draw[fill=white] (.5,.5) circle (.2);
\draw[fill=white] (-.5,.5) circle (.2);
\draw[fill=white] (-.5,-.5) circle (.2);
\draw[fill=white] (.5,-.5) circle (.2);
\node at (.5,.5) {$b$}; % northeast
\node at (-.5,.5) {$a$}; % northwest
\node at (-.5,-.5) {$a$}; % southwest
\node at (.5,-.5) {$b$}; % southeast
\node at (0,.7) {}; % north
\node at (-.7,0) {$k$}; % west
\node at (0,-.7) {}; % south
\node at (.7,0) {$k$}; % east
\end{tikzpicture}
\begin{tikzpicture}
\draw(.5,.5) -- (-.5,.5) -- (-.5,-.5) -- (.5,-.5) -- cycle;
\draw[fill=white] (.5,.5) circle (.2);
\draw[fill=white] (-.5,.5) circle (.2);
\draw[fill=white] (-.5,-.5) circle (.2);
\draw[fill=white] (.5,-.5) circle (.2);
\node at (.5,.5) {$b$}; % northeast
\node at (-.5,.5) {$b$}; % northwest
\node at (-.5,-.5) {$a$}; % southwest
\node at (.5,-.5) {$a$}; % southeast
\node at (0,.7) {$k$}; % north
\node at (-.7,0) {}; % west
\node at (0,-.7) {$k$}; % south
\node at (.7,0) {}; % east
\end{tikzpicture}
\begin{tikzpicture}
\draw(.5,.5) -- (-.5,.5) -- (-.5,-.5) -- (.5,-.5) -- cycle;
\draw[fill=white] (.5,.5) circle (.2);
\draw[fill=white] (-.5,.5) circle (.2);
\draw[fill=white] (-.5,-.5) circle (.2);
\draw[fill=white] (.5,-.5) circle (.2);
\node at (.5,.5) {$a$}; % northeast
\node at (-.5,.5) {$b$}; % northwest
\node at (-.5,-.5) {$b$}; % southwest
\node at (.5,-.5) {$a$}; % southeast
\node at (0,.7) {}; % north
\node at (-.7,0) {$k$}; % west
\node at (0,-.7) {}; % south
\node at (.7,0) {$k$}; % east
\end{tikzpicture}
\vskip6pt%%%%%%%%%%%%%%%%%%%%%%%%%%%%%%%%%%%%%%%%%
\begin{tikzpicture}
\draw(.5,.5) -- (-.5,.5) -- (-.5,-.5) -- (.5,-.5) -- cycle;
\draw[fill=white] (.5,.5) circle (.2);
\draw[fill=white] (-.5,.5) circle (.2);
\draw[fill=white] (-.5,-.5) circle (.2);
\draw[fill=white] (.5,-.5) circle (.2);
\node at (.5,.5) {$a$}; % northeast
\node at (-.5,.5) {$a$}; % northwest
\node at (-.5,-.5) {$b$}; % southwest
\node at (.5,-.5) {$b$}; % southeast
\node at (0,.7) {}; % north
\node at (-.7,0) {$j$}; % west
\node at (0,-.7) {}; % south
\node at (.7,0) {$i$}; % east
\end{tikzpicture}
\begin{tikzpicture}
\draw(.5,.5) -- (-.5,.5) -- (-.5,-.5) -- (.5,-.5) -- cycle;
\draw[fill=white] (.5,.5) circle (.2);
\draw[fill=white] (-.5,.5) circle (.2);
\draw[fill=white] (-.5,-.5) circle (.2);
\draw[fill=white] (.5,-.5) circle (.2);
\node at (.5,.5) {$b$}; % northeast
\node at (-.5,.5) {$a$}; % northwest
\node at (-.5,-.5) {$a$}; % southwest
\node at (.5,-.5) {$b$}; % southeast
\node at (0,.7) {$i$}; % north
\node at (-.7,0) {}; % west
\node at (0,-.7) {$j$}; % south
\node at (.7,0) {}; % east
\end{tikzpicture}
\begin{tikzpicture}
\draw(.5,.5) -- (-.5,.5) -- (-.5,-.5) -- (.5,-.5) -- cycle;
\draw[fill=white] (.5,.5) circle (.2);
\draw[fill=white] (-.5,.5) circle (.2);
\draw[fill=white] (-.5,-.5) circle (.2);
\draw[fill=white] (.5,-.5) circle (.2);
\node at (.5,.5) {$b$}; % northeast
\node at (-.5,.5) {$b$}; % northwest
\node at (-.5,-.5) {$a$}; % southwest
\node at (.5,-.5) {$a$}; % southeast
\node at (0,.7) {}; % north
\node at (-.7,0) {$i$}; % west
\node at (0,-.7) {}; % south
\node at (.7,0) {$j$}; % east
\end{tikzpicture}
\begin{tikzpicture}
\draw(.5,.5) -- (-.5,.5) -- (-.5,-.5) -- (.5,-.5) -- cycle;
\draw[fill=white] (.5,.5) circle (.2);
\draw[fill=white] (-.5,.5) circle (.2);
\draw[fill=white] (-.5,-.5) circle (.2);
\draw[fill=white] (.5,-.5) circle (.2);
\node at (.5,.5) {$a$}; % northeast
\node at (-.5,.5) {$b$}; % northwest
\node at (-.5,-.5) {$b$}; % southwest
\node at (.5,-.5) {$a$}; % southeast
\node at (0,.7) {$j$}; % north
\node at (-.7,0) {}; % west
\node at (0,-.7) {$i$}; % south
\node at (.7,0) {}; % east
\end{tikzpicture}
\caption{Consecutive-letter necklaces not counted in the first collection: $i \not= j$}\label{figure:q2-diagonal-negation}
\end{center}
\end{figure}
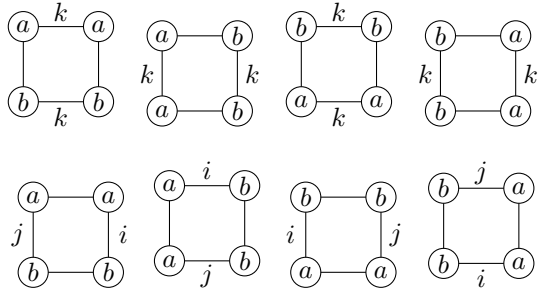
    
\item In the remaining necklaces, either the vertex labels alternate in $a$'s and $b$'s or they do not.

If the $a$'s and $b$'s alternate, using the notation established in Figure~\ref{figure:generic-m4r2}, $s_2=s_4$ so to belong to the second collection according to its original description, $k_1=k_3$ must hold. The necklaces in the second collection where the $a$'s and $b$'s alternate are succinctly depicted in Figure~\ref{figure:q1-north-equals-south}.
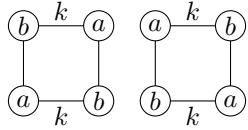
\begin{figure}[hbt]
\begin{center}
\begin{tikzpicture}
\draw(.5,.5) -- (-.5,.5) -- (-.5,-.5) -- (.5,-.5) -- cycle;
\draw[fill=white] (.5,.5) circle (.2);
\draw[fill=white] (-.5,.5) circle (.2);
\draw[fill=white] (-.5,-.5) circle (.2);
\draw[fill=white] (.5,-.5) circle (.2);
\node at (.5,.5) {$a$}; % northeast
\node at (-.5,.5) {$b$}; % northwest
\node at (-.5,-.5) {$a$}; % southwest
\node at (.5,-.5) {$b$}; % southeast
\node at (0,.7) {$k$}; % north
\node at (-.7,0) {}; % west
\node at (0,-.7) {$k$}; % south
\node at (.7,0) {}; % east
\end{tikzpicture}
\begin{tikzpicture}
\draw(.5,.5) -- (-.5,.5) -- (-.5,-.5) -- (.5,-.5) -- cycle;
\draw[fill=white] (.5,.5) circle (.2);
\draw[fill=white] (-.5,.5) circle (.2);
\draw[fill=white] (-.5,-.5) circle (.2);
\draw[fill=white] (.5,-.5) circle (.2);
\node at (.5,.5) {$b$}; % northeast
\node at (-.5,.5) {$a$}; % northwest
\node at (-.5,-.5) {$b$}; % southwest
\node at (.5,-.5) {$a$}; % southeast
\node at (0,.7) {$k$}; % north
\node at (-.7,0) {}; % west
\node at (0,-.7) {$k$}; % south
\node at (.7,0) {}; % east
\end{tikzpicture}
\caption{The $a$'s and $b$'s alternate}\label{figure:q1-north-equals-south}
\end{center}
\end{figure}

The four cases depicted in the first row in Figure~\ref{figure:q2-diagonal-negation} are captured by the condition $k_1=k_3$, or $k_2=k_0$ and $s_2 \not= s_4$.

All in all, the $(4,2,n)$-necklaces in the second collection are precisely the necklaces described by the first row of Figure~\ref{figure:q2-diagonal-negation} or the necklaces described by Figure~\ref{figure:q1-north-equals-south}. These (and only these) are the necklaces accounted for in $Q_1$. That is, a necklace is in the second collection if and only if it is counted in the matrix $Q_1$. We now specify precisely which entry of $Q_1$ accounts for each necklace in the second collection: 

\begin{enumerate}

\item Suppose $k_1 = k_2 = k_3 = k_0$. The $\binom42$ necklaces with uniform edge labeling are counted by the $6$ appearing as the $(\{k_0\},\{k_0\})$-entry (in the main diagonal) of the upper-left block of $Q_1$.

The $\{k_0\}$-entry in the vector $z_1$ is $a_{\{k_0,k_0\}}b_{\{k_0,k_0\}}$, thus the $(\{k_0\},\{k_0\})$-entry $6$ is multiplied by $(a_{\{k_0,k_0\}}b_{\{k_0,k_0\}})^2$ in~\eqref{equation:big-equation}, the monomial represented by necklaces with all edge labels equal to $k_0$.

\item Suppose $k_1 = k_3 \neq k_2 = k_0$. That is, the there are two distinct edge labels $u$ and $v$ used, with the top and bottom edge labels matching, and the left and right edge labels matching. There are $\binom42$ ways to label the vertices with $u$ as the top and bottom edge labels with the remaining edges labeled $v$ and another $\binom42$ ways to label the vertices with $u$ as the left and right edge labels with the remaining edges labeled $v$. These $2\binom42$ necklaces are counted by the $12$ appearing as the $(\{u,v\},\{u,v\})$-entry (in the main diagonal) of the lower-right block of $Q_1$.

The $\{u,v\}$-entry in the vector $z_1$ is $a_{\{u,v\}}b_{\{u,v\}}$, thus the $(\{u,v\},\{u,v\})$-entry $12$ of $Q_1$ is multiplied by $(a_{\{u,v\}}b_{\{u,v\}})^2$ in~\eqref{equation:big-equation}, the monomial represented by necklaces with two distinct edge labels $u$ and $v$, with the top and bottom edges having the same label, and the left and right edges having the same label.

\item Suppose exactly three edge labels are equal. Then the necklace is accounted for in the upper right or bottom left blocks of $Q_1$.

Consider the first row of Figure~\ref{figure:q2-diagonal-negation} or the necklaces described by Figure~\ref{figure:q1-north-equals-south} with the depicted non-adjacent edge labels of $k$. Fix $\ell \not= k$. Six necklaces are obtained by choosing an unmarked edge to be labeled with $k$ and the other unmarked edge to be labeled with $\ell$, and six more ore obtained by swapping the labels on the unmarked edges. These twelve total necklaces are counted by the $6$ appearing as the $(\{k\},\{k,\ell\})$-entry in the upper-right block of $Q_1$ and the $6$ appearing as the $(\{k,\ell\},\{k\})$-entry in the lower-left block of $Q_1$.

The $\{k\}$-entry in the vector $z_1$ is $a_{\{k,k\}}b_{\{k,k\}}$ and the $\{k,\ell\}$-entry is $a_{\{k,\ell\}}b_{\{k,\ell\}}$, thus the $(\{k\},\{k,\ell\})$- and $(\{k,\ell\},\{k\})$-entries are multiplied by $(a_{\{k,k\}}b_{\{k,k\}})(a_{\{k,\ell\}}b_{\{k,\ell\}})$ in~\eqref{equation:big-equation}, the monomial represented by necklaces in the first row of Figure~\ref{figure:q2-diagonal-negation} or the necklaces described by Figure~\ref{figure:q1-north-equals-south} by adding the edge labels $k$ and $\ell$ once each to the unlabeled edges in the figures.

\item Suppose there are three unique edge labels. Then the necklace is accounted for in the non-diagonal part of the bottom right block of $Q_1$.

Consider the first row of Figure~\ref{figure:q2-diagonal-negation} or the necklaces described by Figure~\ref{figure:q1-north-equals-south} with the depicted non-adjacent edge labels of $k$. Fix $u$ and $v$ so that all three of $k,u,v$ are distinct.

Six necklaces are obtained by choosing an unmarked edge to be labeled with $u$ and the other unmarked edge to be labeled with $v$, and six more ore obtained by swapping the labels on the unmarked edges. These twelve total necklaces are counted by the $6$ appearing as the $(\{k,u\},\{k,v\})$-entry in the lower-right block of $Q_1$ and the $6$ appearing as the $(\{k,v\},\{k,u\})$-entry in the lower-right block of $Q_1$. Since $k \not\in \{u,v\}$, these entries are non-diagonal entries in the lower-right block of $Q_1$.

The $\{k,u\}$-entry is $a_{\{k,u\}}b_{\{k,u\}}$ and the $\{k,v\}$-entry is $a_{\{k,v\}}b_{\{k,v\}}$, thus the $(\{k,u\},\{k,v\})$- and $(\{k,v\},\{k,u\})$-entries are multiplied by $(a_{\{k,u\}}b_{\{k,u\}})(a_{\{k,v\}}b_{\{k,v\}})$ in~\eqref{equation:big-equation}, the monomial represented by necklaces in the first row of Figure~\ref{figure:q2-diagonal-negation} or the necklaces described by Figure~\ref{figure:q1-north-equals-south} by adding the edge labels $u$ and $v$ once each to the unlabeled edges in the figures.

\end{enumerate}

By aggregating, for all $u,v$, one can view the $(\{k,u\},\{k,v\})$-entry of $Q_1$ as counting the monomials of the form $(a_{\{k,u\}}b_{\{k,u\}})(a_{\{k,v\}}b_{\{k,v\}})$, but disaggregating based on the cardinalities of $\{k,u,v\}$ and $\{u,v\}$ were helpful due to the differing entries in the matrix $Q_1$.

What remain at this point are the necklaces where the vertex labels alternate $a$'s and $b$'s and the north and south edge labels are not equal, or the necklaces satisfying the constraints of the bottom row of Figure~\ref{figure:q2-diagonal-negation} with the additional constraint that the unmarked edges in the figure must have different edge labels.

\item Since $s_2 \not= s_4$, this considers the four types of necklaces depicted in the bottom row of Figure~\ref{figure:q2-diagonal-negation} with the additional condition that the top and bottom edge labels are not equal. In summary, the necklaces in the third collection are exactly those depicted in Figure~\ref{figure:q2-remaining-main-block}.
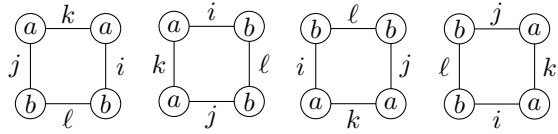
\begin{figure}[hbt]
\begin{center}
\begin{tikzpicture}
\draw(.5,.5) -- (-.5,.5) -- (-.5,-.5) -- (.5,-.5) -- cycle;
\draw[fill=white] (.5,.5) circle (.2);
\draw[fill=white] (-.5,.5) circle (.2);
\draw[fill=white] (-.5,-.5) circle (.2);
\draw[fill=white] (.5,-.5) circle (.2);
\node at (.5,.5) {$a$}; % northeast
\node at (-.5,.5) {$a$}; % northwest
\node at (-.5,-.5) {$b$}; % southwest
\node at (.5,-.5) {$b$}; % southeast
\node at (0,.7) {$k$}; % north
\node at (-.7,0) {$j$}; % west
\node at (0,-.7) {$\ell$}; % south
\node at (.7,0) {$i$}; % east
\end{tikzpicture}
\begin{tikzpicture}
\draw(.5,.5) -- (-.5,.5) -- (-.5,-.5) -- (.5,-.5) -- cycle;
\draw[fill=white] (.5,.5) circle (.2);
\draw[fill=white] (-.5,.5) circle (.2);
\draw[fill=white] (-.5,-.5) circle (.2);
\draw[fill=white] (.5,-.5) circle (.2);
\node at (.5,.5) {$b$}; % northeast
\node at (-.5,.5) {$a$}; % northwest
\node at (-.5,-.5) {$a$}; % southwest
\node at (.5,-.5) {$b$}; % southeast
\node at (0,.7) {$i$}; % north
\node at (-.7,0) {$k$}; % west
\node at (0,-.7) {$j$}; % south
\node at (.7,0) {$\ell$}; % east
\end{tikzpicture}
\begin{tikzpicture}
\draw(.5,.5) -- (-.5,.5) -- (-.5,-.5) -- (.5,-.5) -- cycle;
\draw[fill=white] (.5,.5) circle (.2);
\draw[fill=white] (-.5,.5) circle (.2);
\draw[fill=white] (-.5,-.5) circle (.2);
\draw[fill=white] (.5,-.5) circle (.2);
\node at (.5,.5) {$b$}; % northeast
\node at (-.5,.5) {$b$}; % northwest
\node at (-.5,-.5) {$a$}; % southwest
\node at (.5,-.5) {$a$}; % southeast
\node at (0,.7) {$\ell$}; % north
\node at (-.7,0) {$i$}; % west
\node at (0,-.7) {$k$}; % south
\node at (.7,0) {$j$}; % east
\end{tikzpicture}
\begin{tikzpicture}
\draw(.5,.5) -- (-.5,.5) -- (-.5,-.5) -- (.5,-.5) -- cycle;
\draw[fill=white] (.5,.5) circle (.2);
\draw[fill=white] (-.5,.5) circle (.2);
\draw[fill=white] (-.5,-.5) circle (.2);
\draw[fill=white] (.5,-.5) circle (.2);
\node at (.5,.5) {$a$}; % northeast
\node at (-.5,.5) {$b$}; % northwest
\node at (-.5,-.5) {$b$}; % southwest
\node at (.5,-.5) {$a$}; % southeast
\node at (0,.7) {$j$}; % north
\node at (-.7,0) {$\ell$}; % west
\node at (0,-.7) {$i$}; % south
\node at (.7,0) {$k$}; % east
\end{tikzpicture}
\caption{necklaces in the third collection: $i \not= j$ and $k \not= \ell$}
\label{figure:q2-remaining-main-block}
\end{center}
\end{figure}

These are the necklaces where the vertex labels have adjacent $a$'s and $b$'s, and the conditions $i \neq j$ and $k \neq \ell$ from the figure amount to saying that non-adjacent edge labels are distinct. For a fixed $i,j,k,\ell$ with $i \neq j$ and $k \neq \ell$, there are four necklaces, all of which are obtained as rotations of each other. These necklaces are accounted for by the $4$ appearing as the $(i,j)$-entry in the upper-left or lower-right block of a copy of $Q_2$. 

More specifically, fix $i$ and $j$ such that $i \not= j$. If $k<\ell$ then these necklaces are counted by the $4$ appearing as the $(i,j)$-entry of the upper-left block of the $(k,\ell)$ copy of $Q_2$. If $\ell<k$ then these necklaces are counted by the $4$ appearing as the $(i,j)$-entry of the lower-right block of the $(\ell,k)$ copy of $Q_2$.

If $k < \ell$, the $i$th and $j$th entries in the upper half of the vector $z_{2,(k,\ell)}$ are $a_{\{k,i\}}b_{\{\ell,i\}}$ and $a_{\{k,j\}}b_{\{\ell,j\}}$ respectively, thus the $4$ appearing as the $(i,j)$-entry of the upper-left block of the $(k,\ell)$ copy of $Q_2$ is multiplied by $(a_{\{k,i\}}b_{\{\ell,i\}})(a_{\{k,j\}}b_{\{\ell,j\}})$ in~\eqref{equation:big-equation}, a monomial represented by the four necklaces in Figure~\ref{figure:q2-remaining-main-block}.

If $\ell < k$, the $i$th and $j$th entries in the lower half of the vector $z_{2,(\ell,k)}$ are $a_{\{k,i\}}b_{\{\ell,i\}}$ and $a_{\{k,j\}}b_{\{\ell,j\}}$ respectively, thus the $4$ appearing as the $(i,j)$-entry of the lower-right block of the $(\ell,k)$ copy of $Q_2$ is multiplied by $(a_{\{k,i\}}b_{\{\ell,i\}})(a_{\{k,j\}}b_{\{\ell,j\}})$ in~\eqref{equation:big-equation}, a monomial represented by the four necklaces in Figure~\ref{figure:q2-remaining-main-block}.

The condition $i \not= j$ here implies we have examined non-diagonal entries in the upper-left and lower-right blocks of $Q_2$ matrices. By running through all $i,j,k,\ell$ with $i \neq j$ and $k \neq \ell$, we have used every non-diagonal entry in the upper-left and lower-right blocks of all $\binom{n}{2}$ instances of $Q_2$ in~\eqref{equation:big-equation} exactly once.

\item All necklaces with consecutive $a$'s and $b$'s for vertex labels have been considered. What remains are the necklaces where vertex labels alternate $a$'s and $b$'s. Since the second collection accounted for when the top and bottom edges labels were equal, all that remains are the necklaces with alternating vertex labels with distinct top and bottom edge labels. These are depicted in Figure~\ref{figure:q2-nonmain-block}. These necklaces are counted by all entries in the upper-right and lower-left blocks of all $\binom{n}{2}$ instances of $Q_2$.
\begin{figure}[hbt]
\begin{center}
\begin{tikzpicture}
\draw(.5,.5) -- (-.5,.5) -- (-.5,-.5) -- (.5,-.5) -- cycle;
\draw[fill=white] (.5,.5) circle (.2);
\draw[fill=white] (-.5,.5) circle (.2);
\draw[fill=white] (-.5,-.5) circle (.2);
\draw[fill=white] (.5,-.5) circle (.2);
\node at (.5,.5) {$a$}; % northeast
\node at (-.5,.5) {$b$}; % northwest
\node at (-.5,-.5) {$a$}; % southwest
\node at (.5,-.5) {$b$}; % southeast
\node at (0,.7) {$i$}; % north
\node at (-.7,0) {$\ell$}; % west
\node at (0,-.7) {$j$}; % south
\node at (.7,0) {$k$}; % east
\end{tikzpicture}
\begin{tikzpicture}
\draw(.5,.5) -- (-.5,.5) -- (-.5,-.5) -- (.5,-.5) -- cycle;
\draw[fill=white] (.5,.5) circle (.2);
\draw[fill=white] (-.5,.5) circle (.2);
\draw[fill=white] (-.5,-.5) circle (.2);
\draw[fill=white] (.5,-.5) circle (.2);
\node at (.5,.5) {$b$}; % northeast
\node at (-.5,.5) {$a$}; % northwest
\node at (-.5,-.5) {$b$}; % southwest
\node at (.5,-.5) {$a$}; % southeast
\node at (0,.7) {$i$}; % north
\node at (-.7,0) {$k$}; % west
\node at (0,-.7) {$j$}; % south
\node at (.7,0) {$\ell$}; % east
\end{tikzpicture}
\caption{The $a$'s and $b$'s alternate, with $i \not= j$. There is no condition on the relationship between $k$ and $\ell$.}\label{figure:q2-nonmain-block}
\end{center}
\end{figure}
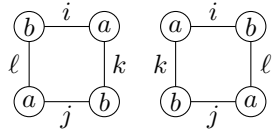

More specifically, if $i < j$ then the two necklaces in Figure~\ref{figure:q2-nonmain-block} are counted by the $2$ appearing as the $(k,\ell)$-entry of the upper-right block of the $(i,j)$ copy of $Q_2$, but if $j < i$ then the two necklaces in Figure~\ref{figure:q2-nonmain-block} are counted by the $2$ appearing as the $(k,\ell)$-entry of the lower-left block of the $(j,i)$ copy of $Q_2$.

If $i < j$, the $k$th entry of the upper half of $z_{2,(i,j)}$ is $a_{\{i,k\}}b_{\{j,k\}}$, and the the $\ell$th entry of the lower half of $z_{2,(i,j)}$ is $a_{\{j,\ell\}}b_{\{i,\ell\}}$, thus the $2$ appearing as the $(k,\ell)$-entry of the upper-right block of the $(i,j)$ copy of $Q_2$ is multiplied by $(a_{\{i,k\}}b_{\{j,k\}})(a_{\{j,\ell\}}b_{\{i,\ell\}})$ in~\eqref{equation:big-equation}, a monomial represented by the two necklaces in Figure~\ref{figure:q2-nonmain-block} satisfying $i < j$.

If $j < i$, the $k$th entry of the lower half of $z_{2,(j,i)}$ is $a_{\{i,k\}}b_{\{j,k\}}$, and the the $\ell$th entry of the upper half of $z_{2,(j,i)}$ is $a_{\{j,\ell\}}b_{\{i,\ell\}}$, thus the $2$ appearing as the $(k,\ell)$-entry of the lower-left block of the $(j,i)$ copy of $Q_2$ is multiplied by $(a_{\{i,k\}}b_{\{j,k\}})(a_{\{j,\ell\}}b_{\{i,\ell\}})$ in~\eqref{equation:big-equation}, a monomial represented by the two necklaces in Figure~\ref{figure:q2-nonmain-block} satisfying $j < i$.

Since the only condition on $i,j,k,\ell$ was that $i \not= j$, we have used every entry in the upper-right and lower-left blocks of all $\binom{n}{2}$ instances of $Q_2$ in~\eqref{equation:big-equation} exactly once.

\end{enumerate}
We partitioned the set of all $(4,2,n)$-necklaces into four so-called collections. Each collection's description was sufficiently manageable that we could enumerate the necklaces therein. In doing so, every non-zero entry of all $1 + \binom{n}{2}$ matrices appearing in~\eqref{equation:big-equation} has been used in the accounting exactly once: the non-zero entries of $Q_1$ counted necklaces in the second collection; the diagonal entries of the $\binom{n}{2}$ instances of $Q_2$ counted necklaces in the first collection, while the non-diagonal entries the upper-left and lower-right blocks of $Q_2$ matrices counted necklaces in the third collection, and the entries in the upper-right and lower-left blocks of $Q_2$ matrices counted necklaces in the fourth collection. Thus, we have proved Theorem~\ref{theorem:main-m4r2}.

\begin{remark}
When thinking about where a necklace is counted in the lower right block of $Q_1$, the indexing of $Q_1$ is helpful. For example, $(\{1,2\},\{2,3\})$-entry of $Q_1$ is $6$, and there are exactly six necklaces with one edge labeled $1$, two edges labeled $2$, and one edge labeled $3$. By writing singletons redundantly, this technique applies to the other blocks of $Q_1$. For example, since $\{1\}=\{1,1\}$, note the $(\{1,1\},\{1,2\})$-entry of $Q_1$ is $6$, and there are exactly six necklaces with three edges labeled $1$ and one edge labeled $2$. This logic was how the second collection was further partitioned based on the block structure of $Q_1$.
\end{remark}

\section{Some remarks}\label{section:m4comments}

Having proved that the coefficient of $t^2$ in $\operatorname{trace}((A+tB)^4)$ is indeed a sum of squares of polynomials in the entries of the $n \times n$ symmetric matrices $A$ and $B$, we discuss some notes of interest.

In finding a sum of squares representation of the coefficient of $t^2$ in $\operatorname{trace}((A+tB)^4)$ for each $n$, we have found a particularly nice non-negative integer solution to an infinite family of semidefinite programs in that the entries of the matrices \emph{count} monomials. Indeed, the sum of the entries in the $1 + \binom{n}{2}$ matrices is
\begin{align*}
6\left[ n + n(n-1) + n(n-1) + \binom{n}{2} \cdot 2(n-1) \right] + \binom{n}{2}[4n^2 + 2n^2 + 2n^2 + 4n^2]
\end{align*}
which is precisely the number $\binom{4}{2}n^4$ of $(4,2,n)$-necklaces, and in turn the number of terms (without collecting) in the coefficient of $t^2$ in $\operatorname{trace}((A+tB)^4)$ as well as the number of monomials appearing in~\eqref{equation:pre-necklace}, again without collecting like terms. Motivated by this, we ask the following question, whose positive answer would be stronger than Conjecture~\ref{conjecture:CDTA}.
\begin{question}\label{question:integer-SOS-with-count}
Given $n \times n$ symmetric matrices $A$ and $B$, if $m \geq r \geq 0$ are even integers, then there exists a positive-semidefinite matrix $Q$ with non-negative integer entries and a vector $z$ of multivariate monomials whose variables are the entries of $A$ and $B$ such that the coefficient of $t^r$ in $\operatorname{trace}((A+tB)^m)$ is $z^TQz$ and the sum of the entries in $Q$ is $\binom{m}{r}n^m$.
\end{question}

\begin{figure}[hbt]
\begin{center}
\begin{tikzpicture}
\draw(.5,.5) -- (-.5,.5) -- (-.5,-.5) -- (.5,-.5) -- cycle;
\draw[fill=white] (.5,.5) circle (.2);
\draw[fill=white] (-.5,.5) circle (.2);
\draw[fill=white] (-.5,-.5) circle (.2);
\draw[fill=white] (.5,-.5) circle (.2);
\node at (.5,.5) {$a$}; % northeast
\node at (-.5,.5) {$a$}; % northwest
\node at (-.5,-.5) {$b$}; % southwest
\node at (.5,-.5) {$b$}; % southeast
\node at (0,.7) {$2$}; % north
\node at (-.7,0) {$2$}; % west
\node at (0,-.7) {$2$}; % south
\node at (.7,0) {$3$}; % east
\end{tikzpicture}
\begin{tikzpicture}
\draw(.5,.5) -- (-.5,.5) -- (-.5,-.5) -- (.5,-.5) -- cycle;
\draw[fill=white] (.5,.5) circle (.2);
\draw[fill=white] (-.5,.5) circle (.2);
\draw[fill=white] (-.5,-.5) circle (.2);
\draw[fill=white] (.5,-.5) circle (.2);
\node at (.5,.5) {$a$}; % northeast
\node at (-.5,.5) {$b$}; % northwest
\node at (-.5,-.5) {$a$}; % southwest
\node at (.5,-.5) {$b$}; % southeast
\node at (0,.7) {$2$}; % north
\node at (-.7,0) {$2$}; % west
\node at (0,-.7) {$3$}; % south
\node at (.7,0) {$2$}; % east
\end{tikzpicture}
\caption{Necklaces in distinct collections associated to the same monomial}\label{figure:monomial-split-q1q2}
\end{center}
\end{figure}
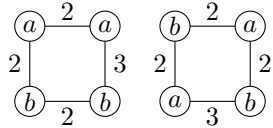

We remark that the proof of Theorem~\ref{theorem:main-m4r2} (and in particular the accounting of necklaces and $Q$-matrix entries in Section~\ref{section:m4-equal-expressions}) is not entirely trivial because a monomial appearing in~\eqref{equation:pre-necklace} may have its accounting ``split'' between various matrices. For example, the two necklaces in Figure~\ref{figure:monomial-split-q1q2} are both associated to $a_{\{2,2\}}b_{\{2,2\}}b_{\{2,3\}}a_{\{2,3\}}$ in~\eqref{equation:pre-necklace} though the former is counted by the $(\{2\},\{2,3\})$-entry in $Q_1$ since it belongs to collection 2 subcollection (c), while the latter is counted by the $(2,2)$-entry in the upper-right block of the $(2,3)$ copy of $Q_2$ since this necklace belongs to the fourth collection.

In the case when $A$ and $B$ are positive semidefinite, then every word in $S_{4,2}(A,B)$ has non-negative trace due to a result of Hillar and Johnson in~\cite{HillarJohnson:n2m5}. This is not true if $A$ and $B$ are only symmetric, even in the case of $2 \times 2$ matrices (even with one diagonal matrix). For instance, if
\[A= \left(\begin{array}{cc}
1 & -3 \\
-3 & 1
\end{array}\right)
\quad\text{and}\quad
B= \left(\begin{array}{cc}
2 & 0 \\
0 & -1
\end{array}\right)
\]
then the trace of $ABAB$ is $-31$, though the trace of $S_{4,2}(A,B)$ is $138$.

Since $A$ and $B$ are both symmetric, by a suitable change of basis, we may have chosen to assume that one of $A$ or $B$ is diagonal. We did not do this so as to obtain in~\eqref{equation:big-equation} an expression that is symmetric in $a$ and $b$ variables which was fruitful in constructing our four-collection classification in Section~\ref{section:m4-equal-expressions}, and moreover, the computational savings obtained were not necessary for us to obtain our SDP feasibility certificate. Of course, we can \emph{ex post facto} present a result for the case when one of $A$ or $B$ is assumed to be diagonal by deleting appropriate rows and columns of $Q_1$ and $Q_2$, and deleting the corresponding entries in the $z$ vectors.

\section{Proof of Theorem~\ref{theorem:main-m8r4}}\label{section:m8r4}

In this section, we prove for every $n \times n$ diagonal matrix $A$ and every $n \times n$ symmetric matrix $B$, if $n \leq 5$, then there exist positive semidefinite matrices $Q_1$, $Q_2$, and $Q_3$ such that the coefficient of $t^4$ in the trace of $(A+tB)^8$ is
\[ z_1^TQ_1z_1 + z_2^TQ_2z_2 + \sum_{1 \leq i < j \leq n} z_{3,(i,j)}^TQ_3z_{3,(i,j)}.\]
The use of diagonal $A$ was needed due to computational limitations which arose in the general case using CVXOPT~\cite{cvxopt} in Sage~\cite{sage}. By suitable simultaneous change of basis, the result applies to general $n \times n$ symmetric matrices where $n \leq 5$.

In~\cite{Collins}, for symmetric $A$ and $B$, Collins, Dykema, and Torres-Ayala prove $S_{2t,2}(A,B)$ and $S_{4t+2,4}(A,B)$ are both cyclically equivalent to a sum of Hermitian squares in the non-commutative variables $A$ and $B$, but $S_{8,4}(A,B)$ is not. We show that the trace of $S_{8,4}(A,B)$ is a sum of squares in the $2n + \binom{n}{2}$ commutative variables $a_{\{i,i\}}$ and $b_{\{i,j\}}$ for $n \leq 5$, which gives progress on the open question of Collins, Dykema, and Torres-Ayala in~\cite{Collins}.

\subsection{Definitions of relevant matrices and monomial vectors}\label{section:m8defineQz}

Given a diagonal matrix $A$ and symmetric matrix $B$, both of size $n \times n$, we denote their entries $a_{\{i,i\}}$ and $b_{\{i,j\}}$.

Let $z_1$ be the vector of size $n$ whose $i$th entry is $a_{\{i,i\}}^2b_{\{i,i\}}^2$. Let $Q_1$ be the $n \times n$ matrix with $70$s on the diagonal, and $0$ elsewhere.

We now describe $z_2$, a vector of size $n(n-1) + \binom{n}{2}$. The $n(n-1)$ entries in the first block of $z_2$ are indexed by ordered pairs $(i,j)$ with distinct entries. The $(i,j)$-entry is $a_{\{i,i\}}^2b_{\{i,j\}}^2$. The $\binom{n}{2}$ entries in the second block of $z_2$ are indexed by $2$-element subsets of $[n]$. The $\{i,j\}$-entry is $a_{\{i,i\}}a_{\{j,j\}}b_{\{i,j\}}^2$. (Note that in the first block, the order matters in the indexing, while in the second block, order is insignificant for indexing. For example, the $(2,3)$-entry in the first block is $a_{\{2,2\}}^2b_{\{2,3\}}^2$ while the $(3,2)$-entry in the first block is $a_{\{3,3\}}^2b_{\{2,3\}}^2$. On the other hand, the $\{2,3\}$-entry of the second block is the $\{3,2\}$-entry, namely $a_{\{2,2\}}a_{\{3,3\}}b_{\{2,3\}}^2$.)

We describe the entries of the $(n(n-1) + \binom{n}{2}) \times (n(n-1) + \binom{n}{2})$ matrix $Q_2$. In the upper left block, diagonal entries are $20$ and all other entries are $0$. In the lower right block, diagonal entries are $36$ and all other entries are $0$. For the non-diagonal blocks, the $(i,j)$-$(\{k,\ell\})$-entry is $16$ if $\{i,j\}=\{k,\ell\}$ and is $0$ otherwise.

For a fixed $1 \leq i < j \leq n$, the vector $z_{3,(i,j)}$ of size $6n-6$ is partitioned into $7$ blocks:
\begin{enumerate}
\item The first block is of size $2$.
Its first entry is $a_{\{i,i\}}^2 b_{\{i,i\}} b_{\{i,j\}}$ and 
its second entry is $a_{\{j,j\}}^2 b_{\{i,i\}} b_{\{j,j\}}$.
\item The second block is of size $2$.
The first entry is $a_{\{i,i\}}a_{\{j,j\}} b_{\{i,i\}} b_{\{i,j\}}$.
The second entry is $a_{\{i,i\}}a_{\{j,j\}} b_{\{j,j\}} b_{\{i,j\}}$.
\item The third block is of size $n-2$.
For each $k \not\in \{i,j\}$, 
there is an entry of the form $a_{\{k,k\}}^2 b_{\{i,k\}} b_{\{j,k\}}$.
\item The fourth block is of size $2(n-2)$.
For each $k \not\in \{i,j\}$, 
there is an entry of the form $a_{\{i,i\}} a_{\{k,k\}} b_{\{i,k\}} b_{\{j,k\}}$
and an entry of the form $a_{\{j,j\}} a_{\{k,k\}} b_{\{i,k\}} b_{\{j,k\}}$.
\item The fifth block is of size $n-1$.
For each $k \not= j$,
there is an entry of the form $a_{\{i,i\}}^2 b_{\{i,k\}} b_{\{j,k\}}$.
While not required for SDP feasibility, the entry where $k=i$ was placed first in this block so that $S_n$ acts on the set of $z_{3,(i,j)}$ vectors.
\item The sixth block is of size $n-2$.
For each $k \not\in \{i,j\}$, 
there is an entry of the form $a_{\{i,i\}} a_{\{j,j\}} b_{\{i,k\}} b_{\{j,k\}}$.
\item The seventh block is of size $n-1$.
For each $k \not= i$,
there is an entry of the form $a_{\{j,j\}}^2 b_{\{i,k\}} b_{\{j,k\}}$.
While not required for SDP feasibility, the entry where $k=j$ was placed first in this block so that $S_n$ acts on the set of $z_{3,(i,j)}$ vectors.
\end{enumerate}

We now describe the $(6n-6) \times (6n-6)$ matrix $Q_3$, which naturally inherits the block structure of the $z_3$ vector. There are $7$ blocks of rows and columns, creating $49$ blocks $(u,v)$ for $u,v = 1,\dots,7$. To provide intuition for how we found $Q_3$, we describe the matrix in blocks using constants and parameters $x_1,\dots,x_{22}$ whose values will be given after the blocks are defined:
\begin{itemize}

\item The $(1,1)$-block of $Q_3$ is the $2 \times 2$ matrix with $120$ on the diagonal and $x_9$ otherwise.

\item The $(1,2)$-block is the $2 \times 2$ matrix with $40$ on the diagonal and $x_1$ otherwise.

\item The $(1,3)$-block is the $2 \times (n-2)$ matrix of all $x_7$.

\item The $(1,4)$-block is the $2 \times 2(n-2)$ matrix with entries $x_{17}$ and $x_{19}$ based on parity, with the upper left entry $x_{17}$. In other words, this block is \[ J_{1 \times (n-2)} \otimes \left(\begin{array}{cc}x_{17} & x_{19} \\ x_{19} & x_{17}\end{array}\right),\] where $J_{1 \times (n-2)}$ is the $1 \times (n-2)$ matrix of all $1$s.

\item The $(1,5)$-block is the $2 \times (n-1)$ matrix with first row all $20$ and second row all $x_4$.
    
\item The $(1,6)$-block is the $2 \times (n-2)$ matrix of all $x_{18}$.

\item The $(1,7)$-block is the $2 \times (n-1)$ matrix with first row all $x_4$ and second row all $20$.

\item The $(2,2)$-block is the $2 \times 2$ matrix with $x_3$ on the diagonal and $x_{10}$ otherwise.

\item The $(2,3)$-block is the $2 \times (n-2)$ matrix of all $x_{11}$.

\item The $(2,4)$-block is the $2 \times 2(n-2)$ matrix with entries $x_{20}$ and $x_{12}$ based on parity, with the upper left entry $x_{20}$.

\item The $(2,5)$-block is the $2 \times (n-1)$ matrix with first row all $x_2$ and second row all $12$.
    
\item The $(2,6)$-block is the $2 \times (n-2)$ matrix of all $x_{8}$.

\item The $(2,7)$-block is the $2 \times (n-1)$ matrix with first row all $12$ and second row all $x_2$.

\item The $(3,3)$-block is the $(n-2) \times (n-2)$ matrix with $40$ on the diagonal and $x_5$ otherwise.

\item The $(3,4)$-block is the $(n-2) \times 2(n-2)$ matrix obtained by taking the $(n-2) \times (n-2)$ matrix with $16$ on the diagonal and $x_{21}$ otherwise and applying the tensor product against $J_{1 \times 2}$.
    
\item The $(3,5)$-block is the $(n-2) \times (n-1)$ matrix of all $4$.

\item The $(3,6)$-block is the $(n-2) \times (n-2)$ matrix with $x_{15}$ on the diagonal and $x_{13}$ otherwise.

\item The $(3,7)$-block is the $(n-2) \times (n-1)$ matrix of all $4$.
    
\item The $(4,4)$-block is the $2(n-2) \times 2(n-2)$ matrix is itself most easily described as a block matrix consisting of $n-2$ blocks of size $2$ each, where the diagonal blocks and the non-diagonal blocks, respectively, consist of the matrices
\[\left(\begin{array}{cc}16 & x_{14} \\ x_{14} & 16\end{array}\right)
\text{ and }
\left(\begin{array}{cc}x_{16} & 2 \\ 2 & x_{16}\end{array}\right),\]
respectively.

\item The $(4,5)$-block is the $2(n-2) \times (n-1)$ matrix whose odd-numbered rows have $8$ and even-numbered rows have $x_{13}$.
    
\item The $(4,6)$-block is the $2(n-2) \times (n-2)$ matrix obtained by taking the $(n-2) \times (n-2)$ matrix with $x_{16}$ on the diagonal and $x_{22}$ otherwise and applying the tensor product againsst $J_{2 \times 1}$.
    
\item The $(4,7)$-block is the $2(n-2) \times (n-1)$ matrix whose odd-numbered rows have $x_{13}$ and even-numbered rows have $8$.
    
\item The $(5,5)$-block is the $(n-1) \times (n-1)$ matrix of all $8$.

\item The $(5,6)$-block is the $(n-1) \times (n-2)$ matrix of all $4$.

\item The $(5,7)$-block is the $(n-1) \times (n-1)$ matrix of all $0$.

\item The $(6,6)$-block is the $(n-2) \times (n-2)$ matrix with $8$ on the diagonal and $x_6$ otherwise.

\item The $(6,7)$-block is the $(n-2) \times (n-1)$ matrix of all $4$.

\item The $(7,7)$-block is the $(n-1) \times (n-1)$ matrix of all $8$.

\end{itemize}
The remaining blocks are completely determined since $Q_3$ is symmetric. For concreteness, when $n=5$, the matrix $Q_3$ described above is
{\tiny
\[
\left(\begin{array}{cc|cc|ccc|cccccc|cccc|ccc|cccc}
120 & x_{9} & 40 & x_{1} & x_{7} & x_{7} & x_{7} & x_{17} & x_{19} & x_{17} & x_{19} & x_{17} & x_{19} & 20 & 20 & 20 & 20 & x_{18} & x_{18} & x_{18} & x_{4} & x_{4} & x_{4} & x_{4} \\
x_{9} & 120 & x_{1} & 40 & x_{7} & x_{7} & x_{7} & x_{19} & x_{17} & x_{19} & x_{17} & x_{19} & x_{17} & x_{4} & x_{4} & x_{4} & x_{4} & x_{18} & x_{18} & x_{18} & 20 & 20 & 20 & 20 \\ \hline
40 & x_{1} & x_{3} & x_{10} & x_{11} & x_{11} & x_{11} & x_{20} & x_{12} & x_{20} & x_{12} & x_{20} & x_{12} & x_{2} & x_{2} & x_{2} & x_{2} & x_{8} & x_{8} & x_{8} & 12 & 12 & 12 & 12 \\
x_{1} & 40 & x_{10} & x_{3} & x_{11} & x_{11} & x_{11} & x_{12} & x_{20} & x_{12} & x_{20} & x_{12} & x_{20} & 12 & 12 & 12 & 12 & x_{8} & x_{8} & x_{8} & x_{2} & x_{2} & x_{2} & x_{2} \\ \hline
x_{7} & x_{7} & x_{11} & x_{11} & 40 & x_{5} & x_{5} & 16 & 16 & x_{21} & x_{21} & x_{21} & x_{21} & 4 & 4 & 4 & 4 & x_{15} & x_{13} & x_{13} & 4 & 4 & 4 & 4 \\
x_{7} & x_{7} & x_{11} & x_{11} & x_{5} & 40 & x_{5} & x_{21} & x_{21} & 16 & 16 & x_{21} & x_{21} & 4 & 4 & 4 & 4 & x_{13} & x_{15} & x_{13} & 4 & 4 & 4 & 4 \\
x_{7} & x_{7} & x_{11} & x_{11} & x_{5} & x_{5} & 40 & x_{21} & x_{21} & x_{21} & x_{21} & 16 & 16 & 4 & 4 & 4 & 4 & x_{13} & x_{13} & x_{15} & 4 & 4 & 4 & 4 \\ \hline
x_{17} & x_{19} & x_{20} & x_{12} & 16 & x_{21} & x_{21} & 16 & x_{14} & x_{16} & 2 & x_{16} & 2 & 8 & 8 & 8 & 8 & x_{16} & x_{22} & x_{22} & x_{13} & x_{13} & x_{13} & x_{13} \\
x_{19} & x_{17} & x_{12} & x_{20} & 16 & x_{21} & x_{21} & x_{14} & 16 & 2 & x_{16} & 2 & x_{16} & x_{13} & x_{13} & x_{13} & x_{13} & x_{16} & x_{22} & x_{22} & 8 & 8 & 8 & 8 \\
x_{17} & x_{19} & x_{20} & x_{12} & x_{21} & 16 & x_{21} & x_{16} & 2 & 16 & x_{14} & x_{16} & 2 & 8 & 8 & 8 & 8 & x_{22} & x_{16} & x_{22} & x_{13} & x_{13} & x_{13} & x_{13} \\
x_{19} & x_{17} & x_{12} & x_{20} & x_{21} & 16 & x_{21} & 2 & x_{16} & x_{14} & 16 & 2 & x_{16} & x_{13} & x_{13} & x_{13} & x_{13} & x_{22} & x_{16} & x_{22} & 8 & 8 & 8 & 8 \\
x_{17} & x_{19} & x_{20} & x_{12} & x_{21} & x_{21} & 16 & x_{16} & 2 & x_{16} & 2 & 16 & x_{14} & 8 & 8 & 8 & 8 & x_{22} & x_{22} & x_{16} & x_{13} & x_{13} & x_{13} & x_{13} \\
x_{19} & x_{17} & x_{12} & x_{20} & x_{21} & x_{21} & 16 & 2 & x_{16} & 2 & x_{16} & x_{14} & 16 & x_{13} & x_{13} & x_{13} & x_{13} & x_{22} & x_{22} & x_{16} & 8 & 8 & 8 & 8 \\ \hline
20 & x_{4} & x_{2} & 12 & 4 & 4 & 4 & 8 & x_{13} & 8 & x_{13} & 8 & x_{13} & 8 & 8 & 8 & 8 & 4 & 4 & 4 & 0 & 0 & 0 & 0 \\
20 & x_{4} & x_{2} & 12 & 4 & 4 & 4 & 8 & x_{13} & 8 & x_{13} & 8 & x_{13} & 8 & 8 & 8 & 8 & 4 & 4 & 4 & 0 & 0 & 0 & 0 \\
20 & x_{4} & x_{2} & 12 & 4 & 4 & 4 & 8 & x_{13} & 8 & x_{13} & 8 & x_{13} & 8 & 8 & 8 & 8 & 4 & 4 & 4 & 0 & 0 & 0 & 0 \\
20 & x_{4} & x_{2} & 12 & 4 & 4 & 4 & 8 & x_{13} & 8 & x_{13} & 8 & x_{13} & 8 & 8 & 8 & 8 & 4 & 4 & 4 & 0 & 0 & 0 & 0 \\ \hline
x_{18} & x_{18} & x_{8} & x_{8} & x_{15} & x_{13} & x_{13} & x_{16} & x_{16} & x_{22} & x_{22} & x_{22} & x_{22} & 4 & 4 & 4 & 4 & 8 & x_{6} & x_{6} & 4 & 4 & 4 & 4 \\
x_{18} & x_{18} & x_{8} & x_{8} & x_{13} & x_{15} & x_{13} & x_{22} & x_{22} & x_{16} & x_{16} & x_{22} & x_{22} & 4 & 4 & 4 & 4 & x_{6} & 8 & x_{6} & 4 & 4 & 4 & 4 \\
x_{18} & x_{18} & x_{8} & x_{8} & x_{13} & x_{13} & x_{15} & x_{22} & x_{22} & x_{22} & x_{22} & x_{16} & x_{16} & 4 & 4 & 4 & 4 & x_{6} & x_{6} & 8 & 4 & 4 & 4 & 4 \\ \hline
x_{4} & 20 & 12 & x_{2} & 4 & 4 & 4 & x_{13} & 8 & x_{13} & 8 & x_{13} & 8 & 0 & 0 & 0 & 0 & 4 & 4 & 4 & 8 & 8 & 8 & 8 \\
x_{4} & 20 & 12 & x_{2} & 4 & 4 & 4 & x_{13} & 8 & x_{13} & 8 & x_{13} & 8 & 0 & 0 & 0 & 0 & 4 & 4 & 4 & 8 & 8 & 8 & 8 \\
x_{4} & 20 & 12 & x_{2} & 4 & 4 & 4 & x_{13} & 8 & x_{13} & 8 & x_{13} & 8 & 0 & 0 & 0 & 0 & 4 & 4 & 4 & 8 & 8 & 8 & 8 \\
x_{4} & 20 & 12 & x_{2} & 4 & 4 & 4 & x_{13} & 8 & x_{13} & 8 & x_{13} & 8 & 0 & 0 & 0 & 0 & 4 & 4 & 4 & 8 & 8 & 8 & 8
\end{array}\right)
\]
}
We postpone for the moment the discussion of why $Q_3$ was presented using parameters $x_1, \dots, x_{22}$. For the specific matrix $Q_3$ we will need for Theorem~\ref{theorem:main-m8r4}, we will set
$x_{1} = 30,
x_{2} = 2,
x_{3} = 40,
x_{4} = 4,
x_{5} = 4,
x_{6} = 4,
x_{7} = 12,
x_{8} = 8,
x_{9} = 24,
x_{10} = 12,
x_{11} = 2,
x_{12} = 12,
x_{13} = 2,
x_{14} = 8,
x_{15} = 4,
x_{16} = 6,
x_{17} = 20,
x_{18} = 10,
x_{19} = 8,
x_{20} = 4,
x_{21} = 2,
x_{22} = 4$.

\subsection{Examples of relevant matrices and monomial vectors}\label{section:m8exampleQz}

We present concrete descriptions of the definitions above when $n=5$. Though these examples are large, we avoid confusion since $n-2 > 2$.

Let $n=5$. Then the transpose of $z_1$ is $(a_{\{1,1\}}^2b_{\{1,1\}}^2, a_{\{2,2\}}^2b_{\{2,2\}}^2, a_{\{3,3\}}^2b_{\{3,3\}}^2, a_{\{4,4\}}^2b_{\{4,4\}}^2, a_{\{5,5\}}^2b_{\{5,5\}}^2)$ and $Q_1$ is the matrix
\[
\left(\begin{array}{rrrrr}
70 & 0 & 0 & 0 & 0 \\
0 & 70 & 0 & 0 & 0 \\
0 & 0 & 70 & 0 & 0 \\
0 & 0 & 0 & 70 & 0 \\
0 & 0 & 0 & 0 & 70
\end{array}\right).
\]
The transpose of the vector $z_2$ is
$(a_{11}^{2} b_{12}^{2},\allowbreak  
a_{11}^{2} b_{13}^{2},\allowbreak  
a_{11}^{2} b_{14}^{2},\allowbreak  
a_{11}^{2} b_{15}^{2},\allowbreak  
a_{22}^{2} b_{12}^{2},\allowbreak  
a_{22}^{2} b_{23}^{2},\allowbreak  
a_{22}^{2} b_{24}^{2},\allowbreak  
a_{22}^{2} b_{25}^{2},\allowbreak  
a_{33}^{2} b_{13}^{2},\allowbreak  
a_{33}^{2} b_{23}^{2},\allowbreak  
a_{33}^{2} b_{34}^{2},\allowbreak  
a_{33}^{2} b_{35}^{2},\allowbreak  
a_{44}^{2} b_{14}^{2},\allowbreak  
a_{44}^{2} b_{24}^{2},\allowbreak  
a_{44}^{2} b_{34}^{2},\allowbreak  
a_{44}^{2} b_{45}^{2},\allowbreak  
a_{55}^{2} b_{15}^{2},\allowbreak  
a_{55}^{2} b_{25}^{2},\allowbreak  
a_{55}^{2} b_{35}^{2},\allowbreak  
a_{55}^{2} b_{45}^{2},\allowbreak   
a_{11} a_{22} b_{12}^{2},\allowbreak  
a_{11} a_{33} b_{13}^{2},\allowbreak  
a_{11} a_{44} b_{14}^{2},\allowbreak  
a_{11} a_{55} b_{15}^{2},\allowbreak  
a_{22} a_{33} b_{23}^{2},\allowbreak  
a_{22} a_{44} b_{24}^{2},\allowbreak  
a_{22} a_{55} b_{25}^{2},\allowbreak  
a_{33} a_{44} b_{34}^{2},\allowbreak  
a_{33} a_{55} b_{35}^{2},\allowbreak  
a_{44} a_{55} b_{45}^{2})$
and the matrix $Q_2$ is
{\tiny
\[
\left(\begin{array}{cccccccccccccccccccc|cccccccccc}
20 & 0 & 0 & 0 & 0 & 0 & 0 & 0 & 0 & 0 & 0 & 0 & 0 & 0 & 0 & 0 & 0 & 0 & 0 & 0 & 16 & 0 & 0 & 0 & 0 & 0 & 0 & 0 & 0 & 0 \\
0 & 20 & 0 & 0 & 0 & 0 & 0 & 0 & 0 & 0 & 0 & 0 & 0 & 0 & 0 & 0 & 0 & 0 & 0 & 0 & 0 & 16 & 0 & 0 & 0 & 0 & 0 & 0 & 0 & 0 \\
0 & 0 & 20 & 0 & 0 & 0 & 0 & 0 & 0 & 0 & 0 & 0 & 0 & 0 & 0 & 0 & 0 & 0 & 0 & 0 & 0 & 0 & 16 & 0 & 0 & 0 & 0 & 0 & 0 & 0 \\
0 & 0 & 0 & 20 & 0 & 0 & 0 & 0 & 0 & 0 & 0 & 0 & 0 & 0 & 0 & 0 & 0 & 0 & 0 & 0 & 0 & 0 & 0 & 16 & 0 & 0 & 0 & 0 & 0 & 0 \\
0 & 0 & 0 & 0 & 20 & 0 & 0 & 0 & 0 & 0 & 0 & 0 & 0 & 0 & 0 & 0 & 0 & 0 & 0 & 0 & 16 & 0 & 0 & 0 & 0 & 0 & 0 & 0 & 0 & 0 \\
0 & 0 & 0 & 0 & 0 & 20 & 0 & 0 & 0 & 0 & 0 & 0 & 0 & 0 & 0 & 0 & 0 & 0 & 0 & 0 & 0 & 0 & 0 & 0 & 16 & 0 & 0 & 0 & 0 & 0 \\
0 & 0 & 0 & 0 & 0 & 0 & 20 & 0 & 0 & 0 & 0 & 0 & 0 & 0 & 0 & 0 & 0 & 0 & 0 & 0 & 0 & 0 & 0 & 0 & 0 & 16 & 0 & 0 & 0 & 0 \\
0 & 0 & 0 & 0 & 0 & 0 & 0 & 20 & 0 & 0 & 0 & 0 & 0 & 0 & 0 & 0 & 0 & 0 & 0 & 0 & 0 & 0 & 0 & 0 & 0 & 0 & 16 & 0 & 0 & 0 \\
0 & 0 & 0 & 0 & 0 & 0 & 0 & 0 & 20 & 0 & 0 & 0 & 0 & 0 & 0 & 0 & 0 & 0 & 0 & 0 & 0 & 16 & 0 & 0 & 0 & 0 & 0 & 0 & 0 & 0 \\
0 & 0 & 0 & 0 & 0 & 0 & 0 & 0 & 0 & 20 & 0 & 0 & 0 & 0 & 0 & 0 & 0 & 0 & 0 & 0 & 0 & 0 & 0 & 0 & 16 & 0 & 0 & 0 & 0 & 0 \\
0 & 0 & 0 & 0 & 0 & 0 & 0 & 0 & 0 & 0 & 20 & 0 & 0 & 0 & 0 & 0 & 0 & 0 & 0 & 0 & 0 & 0 & 0 & 0 & 0 & 0 & 0 & 16 & 0 & 0 \\
0 & 0 & 0 & 0 & 0 & 0 & 0 & 0 & 0 & 0 & 0 & 20 & 0 & 0 & 0 & 0 & 0 & 0 & 0 & 0 & 0 & 0 & 0 & 0 & 0 & 0 & 0 & 0 & 16 & 0 \\
0 & 0 & 0 & 0 & 0 & 0 & 0 & 0 & 0 & 0 & 0 & 0 & 20 & 0 & 0 & 0 & 0 & 0 & 0 & 0 & 0 & 0 & 16 & 0 & 0 & 0 & 0 & 0 & 0 & 0 \\
0 & 0 & 0 & 0 & 0 & 0 & 0 & 0 & 0 & 0 & 0 & 0 & 0 & 20 & 0 & 0 & 0 & 0 & 0 & 0 & 0 & 0 & 0 & 0 & 0 & 16 & 0 & 0 & 0 & 0 \\
0 & 0 & 0 & 0 & 0 & 0 & 0 & 0 & 0 & 0 & 0 & 0 & 0 & 0 & 20 & 0 & 0 & 0 & 0 & 0 & 0 & 0 & 0 & 0 & 0 & 0 & 0 & 16 & 0 & 0 \\
0 & 0 & 0 & 0 & 0 & 0 & 0 & 0 & 0 & 0 & 0 & 0 & 0 & 0 & 0 & 20 & 0 & 0 & 0 & 0 & 0 & 0 & 0 & 0 & 0 & 0 & 0 & 0 & 0 & 16 \\
0 & 0 & 0 & 0 & 0 & 0 & 0 & 0 & 0 & 0 & 0 & 0 & 0 & 0 & 0 & 0 & 20 & 0 & 0 & 0 & 0 & 0 & 0 & 16 & 0 & 0 & 0 & 0 & 0 & 0 \\
0 & 0 & 0 & 0 & 0 & 0 & 0 & 0 & 0 & 0 & 0 & 0 & 0 & 0 & 0 & 0 & 0 & 20 & 0 & 0 & 0 & 0 & 0 & 0 & 0 & 0 & 16 & 0 & 0 & 0 \\
0 & 0 & 0 & 0 & 0 & 0 & 0 & 0 & 0 & 0 & 0 & 0 & 0 & 0 & 0 & 0 & 0 & 0 & 20 & 0 & 0 & 0 & 0 & 0 & 0 & 0 & 0 & 0 & 16 & 0 \\
0 & 0 & 0 & 0 & 0 & 0 & 0 & 0 & 0 & 0 & 0 & 0 & 0 & 0 & 0 & 0 & 0 & 0 & 0 & 20 & 0 & 0 & 0 & 0 & 0 & 0 & 0 & 0 & 0 & 16 \\\hline
16 & 0 & 0 & 0 & 16 & 0 & 0 & 0 & 0 & 0 & 0 & 0 & 0 & 0 & 0 & 0 & 0 & 0 & 0 & 0 & 36 & 0 & 0 & 0 & 0 & 0 & 0 & 0 & 0 & 0 \\
0 & 16 & 0 & 0 & 0 & 0 & 0 & 0 & 16 & 0 & 0 & 0 & 0 & 0 & 0 & 0 & 0 & 0 & 0 & 0 & 0 & 36 & 0 & 0 & 0 & 0 & 0 & 0 & 0 & 0 \\
0 & 0 & 16 & 0 & 0 & 0 & 0 & 0 & 0 & 0 & 0 & 0 & 16 & 0 & 0 & 0 & 0 & 0 & 0 & 0 & 0 & 0 & 36 & 0 & 0 & 0 & 0 & 0 & 0 & 0 \\
0 & 0 & 0 & 16 & 0 & 0 & 0 & 0 & 0 & 0 & 0 & 0 & 0 & 0 & 0 & 0 & 16 & 0 & 0 & 0 & 0 & 0 & 0 & 36 & 0 & 0 & 0 & 0 & 0 & 0 \\
0 & 0 & 0 & 0 & 0 & 16 & 0 & 0 & 0 & 16 & 0 & 0 & 0 & 0 & 0 & 0 & 0 & 0 & 0 & 0 & 0 & 0 & 0 & 0 & 36 & 0 & 0 & 0 & 0 & 0 \\
0 & 0 & 0 & 0 & 0 & 0 & 16 & 0 & 0 & 0 & 0 & 0 & 0 & 16 & 0 & 0 & 0 & 0 & 0 & 0 & 0 & 0 & 0 & 0 & 0 & 36 & 0 & 0 & 0 & 0 \\
0 & 0 & 0 & 0 & 0 & 0 & 0 & 16 & 0 & 0 & 0 & 0 & 0 & 0 & 0 & 0 & 0 & 16 & 0 & 0 & 0 & 0 & 0 & 0 & 0 & 0 & 36 & 0 & 0 & 0 \\
0 & 0 & 0 & 0 & 0 & 0 & 0 & 0 & 0 & 0 & 16 & 0 & 0 & 0 & 16 & 0 & 0 & 0 & 0 & 0 & 0 & 0 & 0 & 0 & 0 & 0 & 0 & 36 & 0 & 0 \\
0 & 0 & 0 & 0 & 0 & 0 & 0 & 0 & 0 & 0 & 0 & 16 & 0 & 0 & 0 & 0 & 0 & 0 & 16 & 0 & 0 & 0 & 0 & 0 & 0 & 0 & 0 & 0 & 36 & 0 \\
0 & 0 & 0 & 0 & 0 & 0 & 0 & 0 & 0 & 0 & 0 & 0 & 0 & 0 & 0 & 16 & 0 & 0 & 0 & 16 & 0 & 0 & 0 & 0 & 0 & 0 & 0 & 0 & 0 & 36
\end{array}\right).
\]}
There are $\binom{5}{2}$ vectors of the form $z_{3,(i,j)}$. Their transposes are
\begin{itemize}
\item $z_{3,(1,2)}^T = (a_{11}^{2} b_{11} b_{12},\allowbreak a_{22}^{2} b_{12} b_{22},\allowbreak a_{11} a_{22} b_{11} b_{12},\allowbreak a_{11} a_{22} b_{12} b_{22},\allowbreak a_{33}^{2} b_{13} b_{23},\allowbreak a_{44}^{2} b_{14} b_{24},\allowbreak a_{55}^{2} b_{15} b_{25},\allowbreak a_{11} a_{33} b_{13} b_{23},\allowbreak a_{22} a_{33} b_{13} b_{23},\allowbreak a_{11} a_{44} b_{14} b_{24},\allowbreak a_{22} a_{44} b_{14} b_{24},\allowbreak a_{11} a_{55} b_{15} b_{25},\allowbreak a_{22} a_{55} b_{15} b_{25},\allowbreak a_{11}^{2} b_{12} b_{22},\allowbreak a_{11}^{2} b_{13} b_{23},\allowbreak a_{11}^{2} b_{14} b_{24},\allowbreak a_{11}^{2} b_{15} b_{25},\allowbreak a_{11} a_{22} b_{13} b_{23},\allowbreak a_{11} a_{22} b_{14} b_{24},\allowbreak a_{11} a_{22} b_{15} b_{25},\allowbreak a_{22}^{2} b_{11} b_{12},\allowbreak a_{22}^{2} b_{13} b_{23},\allowbreak a_{22}^{2} b_{14} b_{24},\allowbreak a_{22}^{2} b_{15} b_{25})$
\item $z_{3,(1,3)}^T = (a_{11}^{2} b_{11} b_{13},\allowbreak a_{33}^{2} b_{13} b_{33},\allowbreak a_{11} a_{33} b_{11} b_{13},\allowbreak a_{11} a_{33} b_{13} b_{33},\allowbreak a_{22}^{2} b_{12} b_{23},\allowbreak a_{44}^{2} b_{14} b_{34},\allowbreak a_{55}^{2} b_{15} b_{35},\allowbreak a_{11} a_{22} b_{12} b_{23},\allowbreak a_{22} a_{33} b_{12} b_{23},\allowbreak a_{11} a_{44} b_{14} b_{34},\allowbreak a_{33} a_{44} b_{14} b_{34},\allowbreak a_{11} a_{55} b_{15} b_{35},\allowbreak a_{33} a_{55} b_{15} b_{35},\allowbreak a_{11}^{2} b_{13} b_{33},\allowbreak a_{11}^{2} b_{12} b_{23},\allowbreak a_{11}^{2} b_{14} b_{34},\allowbreak a_{11}^{2} b_{15} b_{35},\allowbreak a_{11} a_{33} b_{12} b_{23},\allowbreak a_{11} a_{33} b_{14} b_{34},\allowbreak a_{11} a_{33} b_{15} b_{35},\allowbreak a_{33}^{2} b_{11} b_{13},\allowbreak a_{33}^{2} b_{12} b_{23},\allowbreak a_{33}^{2} b_{14} b_{34},\allowbreak a_{33}^{2} b_{15} b_{35})$
\item $z_{3,(1,4)}^T = (a_{11}^{2} b_{11} b_{14},\allowbreak a_{44}^{2} b_{14} b_{44},\allowbreak a_{11} a_{44} b_{11} b_{14},\allowbreak a_{11} a_{44} b_{14} b_{44},\allowbreak a_{22}^{2} b_{12} b_{24},\allowbreak a_{33}^{2} b_{13} b_{34},\allowbreak a_{55}^{2} b_{15} b_{45},\allowbreak a_{11} a_{22} b_{12} b_{24},\allowbreak a_{22} a_{44} b_{12} b_{24},\allowbreak a_{11} a_{33} b_{13} b_{34},\allowbreak a_{33} a_{44} b_{13} b_{34},\allowbreak a_{11} a_{55} b_{15} b_{45},\allowbreak a_{44} a_{55} b_{15} b_{45},\allowbreak a_{11}^{2} b_{14} b_{44},\allowbreak a_{11}^{2} b_{12} b_{24},\allowbreak a_{11}^{2} b_{13} b_{34},\allowbreak a_{11}^{2} b_{15} b_{45},\allowbreak a_{11} a_{44} b_{12} b_{24},\allowbreak a_{11} a_{44} b_{13} b_{34},\allowbreak a_{11} a_{44} b_{15} b_{45},\allowbreak a_{44}^{2} b_{11} b_{14},\allowbreak a_{44}^{2} b_{12} b_{24},\allowbreak a_{44}^{2} b_{13} b_{34},\allowbreak a_{44}^{2} b_{15} b_{45})$
\item $z_{3,(1,5)}^T = (a_{11}^{2} b_{11} b_{15},\allowbreak a_{55}^{2} b_{15} b_{55},\allowbreak a_{11} a_{55} b_{11} b_{15},\allowbreak a_{11} a_{55} b_{15} b_{55},\allowbreak a_{22}^{2} b_{12} b_{25},\allowbreak a_{33}^{2} b_{13} b_{35},\allowbreak a_{44}^{2} b_{14} b_{45},\allowbreak a_{11} a_{22} b_{12} b_{25},\allowbreak a_{22} a_{55} b_{12} b_{25},\allowbreak a_{11} a_{33} b_{13} b_{35},\allowbreak a_{33} a_{55} b_{13} b_{35},\allowbreak a_{11} a_{44} b_{14} b_{45},\allowbreak a_{44} a_{55} b_{14} b_{45},\allowbreak a_{11}^{2} b_{15} b_{55},\allowbreak a_{11}^{2} b_{12} b_{25},\allowbreak a_{11}^{2} b_{13} b_{35},\allowbreak a_{11}^{2} b_{14} b_{45},\allowbreak a_{11} a_{55} b_{12} b_{25},\allowbreak a_{11} a_{55} b_{13} b_{35},\allowbreak a_{11} a_{55} b_{14} b_{45},\allowbreak a_{55}^{2} b_{11} b_{15},\allowbreak a_{55}^{2} b_{12} b_{25},\allowbreak a_{55}^{2} b_{13} b_{35},\allowbreak a_{55}^{2} b_{14} b_{45})$
\item $z_{3,(2,3)}^T = (a_{22}^{2} b_{22} b_{23},\allowbreak a_{33}^{2} b_{23} b_{33},\allowbreak a_{22} a_{33} b_{22} b_{23},\allowbreak a_{22} a_{33} b_{23} b_{33},\allowbreak a_{11}^{2} b_{12} b_{13},\allowbreak a_{44}^{2} b_{24} b_{34},\allowbreak a_{55}^{2} b_{25} b_{35},\allowbreak a_{11} a_{22} b_{12} b_{13},\allowbreak a_{11} a_{33} b_{12} b_{13},\allowbreak a_{22} a_{44} b_{24} b_{34},\allowbreak a_{33} a_{44} b_{24} b_{34},\allowbreak a_{22} a_{55} b_{25} b_{35},\allowbreak a_{33} a_{55} b_{25} b_{35},\allowbreak a_{22}^{2} b_{23} b_{33},\allowbreak a_{22}^{2} b_{12} b_{13},\allowbreak a_{22}^{2} b_{24} b_{34},\allowbreak a_{22}^{2} b_{25} b_{35},\allowbreak a_{22} a_{33} b_{12} b_{13},\allowbreak a_{22} a_{33} b_{24} b_{34},\allowbreak a_{22} a_{33} b_{25} b_{35},\allowbreak a_{33}^{2} b_{22} b_{23},\allowbreak a_{33}^{2} b_{12} b_{13},\allowbreak a_{33}^{2} b_{24} b_{34},\allowbreak a_{33}^{2} b_{25} b_{35})$
\item $z_{3,(2,4)}^T = (a_{22}^{2} b_{22} b_{24},\allowbreak a_{44}^{2} b_{24} b_{44},\allowbreak a_{22} a_{44} b_{22} b_{24},\allowbreak a_{22} a_{44} b_{24} b_{44},\allowbreak a_{11}^{2} b_{12} b_{14},\allowbreak a_{33}^{2} b_{23} b_{34},\allowbreak a_{55}^{2} b_{25} b_{45},\allowbreak a_{11} a_{22} b_{12} b_{14},\allowbreak a_{11} a_{44} b_{12} b_{14},\allowbreak a_{22} a_{33} b_{23} b_{34},\allowbreak a_{33} a_{44} b_{23} b_{34},\allowbreak a_{22} a_{55} b_{25} b_{45},\allowbreak a_{44} a_{55} b_{25} b_{45},\allowbreak a_{22}^{2} b_{24} b_{44},\allowbreak a_{22}^{2} b_{12} b_{14},\allowbreak a_{22}^{2} b_{23} b_{34},\allowbreak a_{22}^{2} b_{25} b_{45},\allowbreak a_{22} a_{44} b_{12} b_{14},\allowbreak a_{22} a_{44} b_{23} b_{34},\allowbreak a_{22} a_{44} b_{25} b_{45},\allowbreak a_{44}^{2} b_{22} b_{24},\allowbreak a_{44}^{2} b_{12} b_{14},\allowbreak a_{44}^{2} b_{23} b_{34},\allowbreak a_{44}^{2} b_{25} b_{45})$
\item $z_{3,(2,5)}^T = (a_{22}^{2} b_{22} b_{25},\allowbreak a_{55}^{2} b_{25} b_{55},\allowbreak a_{22} a_{55} b_{22} b_{25},\allowbreak a_{22} a_{55} b_{25} b_{55},\allowbreak a_{11}^{2} b_{12} b_{15},\allowbreak a_{33}^{2} b_{23} b_{35},\allowbreak a_{44}^{2} b_{24} b_{45},\allowbreak a_{11} a_{22} b_{12} b_{15},\allowbreak a_{11} a_{55} b_{12} b_{15},\allowbreak a_{22} a_{33} b_{23} b_{35},\allowbreak a_{33} a_{55} b_{23} b_{35},\allowbreak a_{22} a_{44} b_{24} b_{45},\allowbreak a_{44} a_{55} b_{24} b_{45},\allowbreak a_{22}^{2} b_{25} b_{55},\allowbreak a_{22}^{2} b_{12} b_{15},\allowbreak a_{22}^{2} b_{23} b_{35},\allowbreak a_{22}^{2} b_{24} b_{45},\allowbreak a_{22} a_{55} b_{12} b_{15},\allowbreak a_{22} a_{55} b_{23} b_{35},\allowbreak a_{22} a_{55} b_{24} b_{45},\allowbreak a_{55}^{2} b_{22} b_{25},\allowbreak a_{55}^{2} b_{12} b_{15},\allowbreak a_{55}^{2} b_{23} b_{35},\allowbreak a_{55}^{2} b_{24} b_{45})$
\item $z_{3,(3,4)}^T = (a_{33}^{2} b_{33} b_{34},\allowbreak a_{44}^{2} b_{34} b_{44},\allowbreak a_{33} a_{44} b_{33} b_{34},\allowbreak a_{33} a_{44} b_{34} b_{44},\allowbreak a_{11}^{2} b_{13} b_{14},\allowbreak a_{22}^{2} b_{23} b_{24},\allowbreak a_{55}^{2} b_{35} b_{45},\allowbreak a_{11} a_{33} b_{13} b_{14},\allowbreak a_{11} a_{44} b_{13} b_{14},\allowbreak a_{22} a_{33} b_{23} b_{24},\allowbreak a_{22} a_{44} b_{23} b_{24},\allowbreak a_{33} a_{55} b_{35} b_{45},\allowbreak a_{44} a_{55} b_{35} b_{45},\allowbreak a_{33}^{2} b_{34} b_{44},\allowbreak a_{33}^{2} b_{13} b_{14},\allowbreak a_{33}^{2} b_{23} b_{24},\allowbreak a_{33}^{2} b_{35} b_{45},\allowbreak a_{33} a_{44} b_{13} b_{14},\allowbreak a_{33} a_{44} b_{23} b_{24},\allowbreak a_{33} a_{44} b_{35} b_{45},\allowbreak a_{44}^{2} b_{33} b_{34},\allowbreak a_{44}^{2} b_{13} b_{14},\allowbreak a_{44}^{2} b_{23} b_{24},\allowbreak a_{44}^{2} b_{35} b_{45})$
\item $z_{3,(3,5)}^T = (a_{33}^{2} b_{33} b_{35},\allowbreak a_{55}^{2} b_{35} b_{55},\allowbreak a_{33} a_{55} b_{33} b_{35},\allowbreak a_{33} a_{55} b_{35} b_{55},\allowbreak a_{11}^{2} b_{13} b_{15},\allowbreak a_{22}^{2} b_{23} b_{25},\allowbreak a_{44}^{2} b_{34} b_{45},\allowbreak a_{11} a_{33} b_{13} b_{15},\allowbreak a_{11} a_{55} b_{13} b_{15},\allowbreak a_{22} a_{33} b_{23} b_{25},\allowbreak a_{22} a_{55} b_{23} b_{25},\allowbreak a_{33} a_{44} b_{34} b_{45},\allowbreak a_{44} a_{55} b_{34} b_{45},\allowbreak a_{33}^{2} b_{35} b_{55},\allowbreak a_{33}^{2} b_{13} b_{15},\allowbreak a_{33}^{2} b_{23} b_{25},\allowbreak a_{33}^{2} b_{34} b_{45},\allowbreak a_{33} a_{55} b_{13} b_{15},\allowbreak a_{33} a_{55} b_{23} b_{25},\allowbreak a_{33} a_{55} b_{34} b_{45},\allowbreak a_{55}^{2} b_{33} b_{35},\allowbreak a_{55}^{2} b_{13} b_{15},\allowbreak a_{55}^{2} b_{23} b_{25},\allowbreak a_{55}^{2} b_{34} b_{45})$
\item $z_{3,(4,5)}^T = (a_{44}^{2} b_{44} b_{45},\allowbreak a_{55}^{2} b_{45} b_{55},\allowbreak a_{44} a_{55} b_{44} b_{45},\allowbreak a_{44} a_{55} b_{45} b_{55},\allowbreak a_{11}^{2} b_{14} b_{15},\allowbreak a_{22}^{2} b_{24} b_{25},\allowbreak a_{33}^{2} b_{34} b_{35},\allowbreak a_{11} a_{44} b_{14} b_{15},\allowbreak a_{11} a_{55} b_{14} b_{15},\allowbreak a_{22} a_{44} b_{24} b_{25},\allowbreak a_{22} a_{55} b_{24} b_{25},\allowbreak a_{33} a_{44} b_{34} b_{35},\allowbreak a_{33} a_{55} b_{34} b_{35},\allowbreak a_{44}^{2} b_{45} b_{55},\allowbreak a_{44}^{2} b_{14} b_{15},\allowbreak a_{44}^{2} b_{24} b_{25},\allowbreak a_{44}^{2} b_{34} b_{35},\allowbreak a_{44} a_{55} b_{14} b_{15},\allowbreak a_{44} a_{55} b_{24} b_{25},\allowbreak a_{44} a_{55} b_{34} b_{35},\allowbreak a_{55}^{2} b_{44} b_{45},\allowbreak a_{55}^{2} b_{14} b_{15},\allowbreak a_{55}^{2} b_{24} b_{25},\allowbreak a_{55}^{2} b_{34} b_{35})$
\end{itemize}
and based on the specified values of the parameters $x_1, \dots, x_{22}$, the matrix $Q_3$ is
\[
\left(\begin{array}{cc|cc|ccc|cccccc|cccc|ccc|cccc}
120 & 24 & 40 & 30 & 12 & 12 & 12 & 20 & 8 & 20 & 8 & 20 & 8 & 20 & 20 & 20 & 20 & 10 & 10 & 10 & 4 & 4 & 4 & 4 \\
24 & 120 & 30 & 40 & 12 & 12 & 12 & 8 & 20 & 8 & 20 & 8 & 20 & 4 & 4 & 4 & 4 & 10 & 10 & 10 & 20 & 20 & 20 & 20 \\ \hline
40 & 30 & 40 & 12 & 2 & 2 & 2 & 4 & 12 & 4 & 12 & 4 & 12 & 2 & 2 & 2 & 2 & 8 & 8 & 8 & 12 & 12 & 12 & 12 \\
30 & 40 & 12 & 40 & 2 & 2 & 2 & 12 & 4 & 12 & 4 & 12 & 4 & 12 & 12 & 12 & 12 & 8 & 8 & 8 & 2 & 2 & 2 & 2 \\ \hline
12 & 12 & 2 & 2 & 40 & 4 & 4 & 16 & 16 & 2 & 2 & 2 & 2 & 4 & 4 & 4 & 4 & 4 & 2 & 2 & 4 & 4 & 4 & 4 \\
12 & 12 & 2 & 2 & 4 & 40 & 4 & 2 & 2 & 16 & 16 & 2 & 2 & 4 & 4 & 4 & 4 & 2 & 4 & 2 & 4 & 4 & 4 & 4 \\
12 & 12 & 2 & 2 & 4 & 4 & 40 & 2 & 2 & 2 & 2 & 16 & 16 & 4 & 4 & 4 & 4 & 2 & 2 & 4 & 4 & 4 & 4 & 4 \\ \hline
20 & 8 & 4 & 12 & 16 & 2 & 2 & 16 & 8 & 6 & 2 & 6 & 2 & 8 & 8 & 8 & 8 & 6 & 4 & 4 & 2 & 2 & 2 & 2 \\
8 & 20 & 12 & 4 & 16 & 2 & 2 & 8 & 16 & 2 & 6 & 2 & 6 & 2 & 2 & 2 & 2 & 6 & 4 & 4 & 8 & 8 & 8 & 8 \\
20 & 8 & 4 & 12 & 2 & 16 & 2 & 6 & 2 & 16 & 8 & 6 & 2 & 8 & 8 & 8 & 8 & 4 & 6 & 4 & 2 & 2 & 2 & 2 \\
8 & 20 & 12 & 4 & 2 & 16 & 2 & 2 & 6 & 8 & 16 & 2 & 6 & 2 & 2 & 2 & 2 & 4 & 6 & 4 & 8 & 8 & 8 & 8 \\
20 & 8 & 4 & 12 & 2 & 2 & 16 & 6 & 2 & 6 & 2 & 16 & 8 & 8 & 8 & 8 & 8 & 4 & 4 & 6 & 2 & 2 & 2 & 2 \\
8 & 20 & 12 & 4 & 2 & 2 & 16 & 2 & 6 & 2 & 6 & 8 & 16 & 2 & 2 & 2 & 2 & 4 & 4 & 6 & 8 & 8 & 8 & 8 \\\hline
20 & 4 & 2 & 12 & 4 & 4 & 4 & 8 & 2 & 8 & 2 & 8 & 2 & 8 & 8 & 8 & 8 & 4 & 4 & 4 & 0 & 0 & 0 & 0 \\
20 & 4 & 2 & 12 & 4 & 4 & 4 & 8 & 2 & 8 & 2 & 8 & 2 & 8 & 8 & 8 & 8 & 4 & 4 & 4 & 0 & 0 & 0 & 0 \\
20 & 4 & 2 & 12 & 4 & 4 & 4 & 8 & 2 & 8 & 2 & 8 & 2 & 8 & 8 & 8 & 8 & 4 & 4 & 4 & 0 & 0 & 0 & 0 \\
20 & 4 & 2 & 12 & 4 & 4 & 4 & 8 & 2 & 8 & 2 & 8 & 2 & 8 & 8 & 8 & 8 & 4 & 4 & 4 & 0 & 0 & 0 & 0 \\ \hline
10 & 10 & 8 & 8 & 4 & 2 & 2 & 6 & 6 & 4 & 4 & 4 & 4 & 4 & 4 & 4 & 4 & 8 & 4 & 4 & 4 & 4 & 4 & 4 \\
10 & 10 & 8 & 8 & 2 & 4 & 2 & 4 & 4 & 6 & 6 & 4 & 4 & 4 & 4 & 4 & 4 & 4 & 8 & 4 & 4 & 4 & 4 & 4 \\
10 & 10 & 8 & 8 & 2 & 2 & 4 & 4 & 4 & 4 & 4 & 6 & 6 & 4 & 4 & 4 & 4 & 4 & 4 & 8 & 4 & 4 & 4 & 4 \\ \hline
4 & 20 & 12 & 2 & 4 & 4 & 4 & 2 & 8 & 2 & 8 & 2 & 8 & 0 & 0 & 0 & 0 & 4 & 4 & 4 & 8 & 8 & 8 & 8 \\
4 & 20 & 12 & 2 & 4 & 4 & 4 & 2 & 8 & 2 & 8 & 2 & 8 & 0 & 0 & 0 & 0 & 4 & 4 & 4 & 8 & 8 & 8 & 8 \\
4 & 20 & 12 & 2 & 4 & 4 & 4 & 2 & 8 & 2 & 8 & 2 & 8 & 0 & 0 & 0 & 0 & 4 & 4 & 4 & 8 & 8 & 8 & 8 \\
4 & 20 & 12 & 2 & 4 & 4 & 4 & 2 & 8 & 2 & 8 & 2 & 8 & 0 & 0 & 0 & 0 & 4 & 4 & 4 & 8 & 8 & 8 & 8
\end{array}\right).
\]

\subsection{Proofs of positive semidefiniteness}\label{section:m8provePSD}

Since $Q_1$ is a positive scaling of the identity matrix, it is clearly positive semidefinite. The upper left block of $Q_2$ is clearly invertible. Its Schur complement (see, for instance, Section A.5.5 of~\cite{BoydVandenberghe}) is a diagonal matrix where every entry is $36 - \frac1{20}(16^2+16^2) = \frac{52}{5} > 0$, thus $Q_2$ is positive definite. For the case of $n=5$, the characteristic polynomial of $Q_3$ is
\begin{dmath*}
x^{24} - 624x^{23} + 153560x^{22} - 20114720x^{21} + 1579524768x^{20} - 78655770752x^{19} + 2550222470144x^{18} - 54517284561408x^{17} + 777531250595072x^{16} - 7542940488812544x^{15} + 50698374763948032x^{14} - 238763500839682048x^{13} + 788811023799615488x^{12} - 1807444411797995520x^{11} + 2792148062679072768x^{10} - 2752683589605785600x^{9} + 1559961443934142464x^{8} - 408038627117891584x^{7} + 28112973650198528x^{6}.
\end{dmath*}
Thus, when $n=5$, the eigenvalues of $Q_3$ are $0$ (with multiplicity $6$), $0.10251621889932427$, $0.5316579744838135$, $1.772545830208574$, $2.437392431363524$ (with multiplicity $2$), $2.878415693912668$, $4$ (with multiplicity $2$), $5.152819633092394$ (with multiplicity $2$), $7.600606876273283$, $28.20275708078755$, $48.40978793554409$ (with multiplicity $2$), $52.92420941747237$, $53.91800659262472$, $116.7715867778353$, and $239.2976975375025$ with multiplicities. Submatrices of this matrix, obtained by selecting the appropriate [corresponding] rows and columns yield the $Q_3$ matrix for $n < 5$ and are thus positive semidefinite.

\subsection{Proof of relevance to \eqref{equation:big-equation-m8}}\label{section:m8-equal-expressions}

Though Section~\ref{section:m8exampleQz} gave specific values for the parameters $x_1, \dots, x_{22}$, in fact, the interested reader can verify that, for all $n$, the coefficient of $t^4$ in the trace of $(A+tB)^8$ will equal~\eqref{equation:big-equation-m8} if the parameters used to define the matrix $Q_3$ satisfy the following system of linear equations:
\begin{align*}
&x_{1} + x_{2} = 32\\
&x_{3} + 2 x_{4} = 48\\
&x_{5} + x_{6} = 8\\
&x_{4} + x_{7} + x_{8} = 24\\
&x_{9} + x_{10} = 36\\
&x_{11} + x_{12} + x_{13} = 16\\
&x_{14} + x_{15} = 12\\
&x_{13} + x_{16} = 8\\
&x_{2} + x_{17} + x_{18} = 32\\
&x_{19} + x_{20} = 12\\
&x_{13} + x_{21} + x_{22} = 8.
\end{align*}
Via a symbolic algebra system, one can verify that if $Q_3$ is defined according to Section~\ref{section:m8exampleQz} in such a way that the system of linear equations above is satisfied, then the coefficient of $t^4$ in the trace of $(A+tB)^8$ will equal~\eqref{equation:big-equation-m8}. This has been verified computationally up to $n=9$.

\section{Concluding remarks}\label{section:last}

Theorem~\ref{theorem:main-m8r4} provides positive evidence that Conjecture~\ref{conjecture:CDTA} from~\cite{Collins} is true in the still open case of $(m,r)=(8,4)$. Though we considered the case when $A$ is diagonal, our result in the $(m,r)=(8,4)$ case gives evidence that an appropriately-modified variant of Question~\ref{question:integer-SOS-with-count} is true, namely where the sum of the entries is $\binom{m}{r}n^r$. In fact, when $m=8$ and $r=4$, one can prove that for all $n$, the sum of the $Q$-matrix entries is precisely $\binom{m}{r}\sqrt{n^m}$ if the system of linear equations in Section~\ref{section:m8-equal-expressions} is satisfied.

One could ask for a proof of the content in Section~\ref{section:m8-equal-expressions}, perhaps modeled after the proof given in Section~\ref{section:m4-equal-expressions}: account for all $(8,4,n)$-necklaces (assuming $A$ is diagonal, creating restrictions on edge labels) using matrix entries. To do this, the structure of $Q_1$ and $Q_2$ are very clear, but the structure of $Q_3$ is not as clear.

This leads us to now discuss the long-promised purpose behind presenting $Q_3$ with parameters. In the interest of finding a description of $Q_3$ for Theorem~\ref{theorem:main-m8r4} for all $n$, we applied several Ans\"{a}tze based on $S_n$ symmetry. These led to certain entries of $Q_3$ being completely determined, and those entries of $Q_3$ were presented in Section~\ref{section:m8defineQz} as constants. Additionally, the symmetries led us to declaring that certain entries of $Q_3$ should be equal, leading to our $22$ parameters. Without much extra work, through coefficient matching, we constructed the system of linear equations in Section~\ref{section:m8-equal-expressions} as a necessary condition for Theorem~\ref{theorem:main-m8r4} under our Ans\"{a}tze. Since we wanted the entries of our matrices to be non-negative due to Question~\ref{question:integer-SOS-with-count}, this in turn meant we needed all parameters to be non-negative. These constraints together with $Q_3 \succeq 0$ can be easily implemented by a suitable multi-block SDP.

Naturally, it is desirable to consider Conjecture~\ref{conjecture:CDTA} and Question~\ref{question:integer-SOS-with-count} for the cases where $r > 4$ even. While work on cyclic equivalence to a sum of Hermitian squares in non-commutative variables has been fruitful, Theorems~\ref{theorem:main-m4r2} and~\ref{theorem:main-m8r4} ask us to at least consider whether these questions should be revisited using commutative variables. Of course, the sizes of the SDPs grow considerably as $m$, $r$, and $n$ grow. These cases could be considered following~\cite{GP, Vallentin}, though the projection formulae onto invariant subspaces used even in the early algorithms may mean that the integrality matter in Question~\ref{question:integer-SOS-with-count} may be challenging to recover.

The entire exposition here could be stated for complex-valued Hermitian matrices $A$ and $B$ without much extra effort. In every instance of $z^TQz$, one of the $z$ vectors will need to be conjugated. We leave the details to the interested reader.

We thank Frank Vallentin and Frederik von Heymann for helpful conversations in the initial stages of investigation.

\end{document}